%
%
%
%
\hsize=5in
\baselineskip=12pt
\vsize=20.3cm
\parindent=10pt
\pretolerance=40
\predisplaypenalty=0
\displaywidowpenalty=0
\finalhyphendemerits=0
\hfuzz=2pt
\frenchspacing
\footline={\ifnum\pageno=1\else\hfil\tenrm\number\pageno\hfil\fi}
\hyphenation{semi-prime}
\hyphenation{group-like}
\hyphenation{group-likes}
%
%
\input amssym.def
\font\tenbbold=bbold10
\font\sevenbbold=bbold7
\newfam\bbfam
\textfont\bbfam=\tenbbold
\scriptfont\bbfam=\sevenbbold
\def\titlefonts{\baselineskip=1.44\baselineskip
	\font\titlef=cmbx12
	\titlef
	}
\font\tenib=cmmib10 
\skewchar\tenib='177
\def\boldfonts{\bf
	\textfont0=\tenbf
	\textfont1=\tenib
	}
\font\ninerm=cmr9
\font\ninebf=cmbx9
\font\ninei=cmmi9
\skewchar\ninei='177
\font\ninesy=cmsy9
\skewchar\ninesy='60
\font\nineit=cmti9
\def\reffonts{\baselineskip=0.9\baselineskip
	\textfont0=\ninerm
	\def\rm{\fam0\ninerm}%
	\textfont1=\ninei
	\textfont2=\ninesy
	\textfont\bffam=\ninebf
	\def\bf{\fam\bffam\ninebf}%
	\def\it{\nineit}%
	}
%
%
\def\frontmatter{\vbox{}\vskip1cm\bgroup
	\leftskip=0pt plus1fil\rightskip=0pt plus1fil
	\parindent=0pt
	\parfillskip=0pt
	\pretolerance=10000
	}
\def\endfrontmatter{\egroup\bigskip}
\def\title#1{{\titlefonts#1\par}}
\def\author#1{\bigskip#1\par}
\def\address#1{\medskip{\reffonts\it#1}}
\def\email#1{\smallskip{\reffonts{\it E-mail: }\rm#1}}
\def\section#1\par{\ifdim\lastskip<\bigskipamount\removelastskip\fi
	\penalty-250\bigskip
	\vbox{\leftskip=0pt plus1fil\rightskip=0pt plus1fil
	\parindent=0pt
	\parfillskip=0pt
	\pretolerance=10000{\boldfonts#1}}\nobreak\medskip
	}
\def\emph#1{{\it#1}}
\def\proclaim#1. {\medbreak\bgroup{\noindent\bf#1.}\ \it}
\def\endproclaim{\egroup
	\ifdim\lastskip<\medskipamount\removelastskip\medskip\fi}
\newdimen\itemsize
\def\setitemsize#1 {{\setbox0\hbox{#1\ }
	\global\itemsize=\wd0}}
\def\item#1 #2\par{\ifdim\lastskip<\smallskipamount\removelastskip\smallskip\fi
	{\leftskip=\itemsize
	\noindent\hskip-\leftskip
	\hbox to\leftskip{\hfil\rm#1\ }#2\par}\smallskip}
\def\Proof#1. {\ifdim\lastskip<\medskipamount\removelastskip\medskip\fi
	{\noindent\it Proof\if\space#1\space\else\ \fi#1.}\ }
\def\endproof{\hfill\hbox{}\quad\hbox{}\hfill\llap{$\square$}\medskip}
%
%
\newcount\citation
\newtoks\citetoks
\def\citedef#1\endcitedef{\citetoks={#1\endcitedef}}
\def\endcitedef#1\endcitedef{}
\def\citenum#1{\citation=0\def\curcite{#1}%
	\expandafter\checkendcite\the\citetoks}
\def\checkendcite#1{\ifx\endcitedef#1?\else
	\expandafter\lookcite\expandafter#1\fi}
\def\lookcite#1 {\advance\citation by1\def\auxcite{#1}%
	\ifx\auxcite\curcite\the\citation\expandafter\endcitedef\else
	\expandafter\checkendcite\fi}
\def\cite#1{\makecite#1,\cite}
\def\makecite#1,#2{[\citenum{#1}\ifx\cite#2]\else\expandafter\clearcite\expandafter#2\fi}
\def\clearcite#1,\cite{, #1]}
%
%
\def\references{\section References\par
	\bgroup
	\parindent=0pt
	\reffonts
	\rm
	\frenchspacing
	\setbox0\hbox{99. }\leftskip=\wd0
	}
\def\endreferences{\egroup}
\newtoks\authtoks
\newif\iffirstauth
\def\checkendauth#1{\ifx\auth#1%
    \iffirstauth\the\authtoks
    \else{} and \the\authtoks\fi,%
  \else\iffirstauth\the\authtoks\firstauthfalse
    \else, \the\authtoks\fi
    \expandafter\nextauth\expandafter#1\fi
}
\def\nextauth#1,#2;{\authtoks={#1 #2}\checkendauth}
\def\auth#1{\nextauth#1;\auth}
\newif\ifinbook
\newif\ifbookref
\def\nextref#1 {\par\hskip-\leftskip
	\hbox to\leftskip{\hfil\citenum{#1}.\ }%
	\initnextref}
\def\initnextref{\bookreffalse\inbookfalse\firstauthtrue\ignorespaces}
\def\paper#1{{\it#1},}

\def\book#1{\bookreftrue{\it#1},}
\def\journal#1{#1\ifinbook,\fi}
\def\BkSer#1{#1,}
\def\Vol#1{{\bf#1}}
\def\BkVol#1{Vol. #1,}
\def\publisher#1{#1,}
\def\Year#1{\ifbookref #1.\else\ifinbook #1,\else(#1)\fi\fi}
\def\Pages#1{\makepages#1.}
\long\def\makepages#1-#2.#3{\ifinbook pp. \fi#1--#2\ifx\par#3.\fi#3}
\def\inRus{{ \rm(in Russian)}}
\def\etransl#1{English translation in \journal{#1}}
%
%
\newsymbol\square 1003
\let\Rar\Rightarrow
\let\Lrar\Leftrightarrow
\let\rhu\rightharpoonup
\let\rhd\rightharpoondown
\let\ot\otimes
\let\sbs\subset
\let\sps\supset
\newsymbol\varnothing 203F
\let\noth\varnothing
\def\:{\,\mathord:\,}
\def\ann{\mathop{\rm ann}\nolimits}
\def\chr{\mathop{\rm char}\nolimits}
\def\co#1{^{\mkern1mu{\rm co}\mkern1mu#1}}
\def\cop{^{\rm cop}}
\def\End{\mathop{\rm End}\nolimits}
\def\Hom{\mathop{\rm Hom}\nolimits}
\def\id{{\rm id}}
\def\Id{\mathop{\rm Id}\nolimits}
\def\Im{\mathop{\rm Im}}
\def\Ker{\mathop{\rm Ker}}
\def\lng{\mathop{\rm length}\nolimits}
\def\Max{\mathop{\rm Max}\nolimits}
\def\mod#1{\ifinner\mskip8mu(\mathop{\rm mod}#1)
        \else\mskip12mu(\mathop{\rm mod}#1)\fi}
\def\op{^{\rm op}}
\def\rank{\mathop{\rm rank}\nolimits}
\def\Spec{\mathop{\rm Spec}}
\def\tsum{\textstyle\sum\limits}
\def\lmapr#1#2{{}\mathrel{\smash{\mathop{\count0=#1
  \loop
    \ifnum\count0>0
    \advance\count0 by-1\smash{\mathord-}\mkern-4mu
  \repeat
  \mathord\rightarrow}\limits^{#2}}}{}}
\def\lmapd#1#2#3{\llap{$\vcenter{\hbox{$\scriptstyle{#2}$}}$}
  \left\downarrow\vcenter to#1{}\right.
  \rlap{$\vcenter{\hbox{$\scriptstyle{#3}$}}$}}
\def\diagram#1{\vbox{\halign{&\hfil$##$\hfil\cr #1}}}
\let\al\alpha
\let\ga\gamma
\let\de\delta
\let\ep\varepsilon
\let\io\iota
\let\la\lambda
\let\ph\varphi
\let\ze\zeta
\let\De\Delta
\let\Si\Sigma
\def\bbQ{{\Bbb Q}}
\def\hatt{{\hat t}}
\def\qov{{\overline q}}
\def\Aov{{\overline A}}
\def\Qov{{\overline Q}}

\def\bbk{{\fam\bbfam k}}
\def\0{_{(0)}}
\def\1{_{(1)}}
\def\2{_{(2)}}
\def\AHA{A\ot_{A^H}\!A}
\def\AM{\vphantom{\M}_A\M}
\def\Ascr{{\cal A}}
\def\BM{\vphantom{\M}_B\M}
\def\Bscr{{\cal B}}
\def\E{{\cal E}}
\def\Endk{\End_{\,\bbk}}
\def\F{{\cal F}}
\def\FHA{\F_H(A)}
\def\I{{\cal I}}
\def\HAM{H\hbox{-}\vphantom{\M}_A\M}
\def\HAMA{\HAM_A}
\def\HBM{H\hbox{-}\vphantom{\M}_B\M}
\def\HMA{H\hbox{-}\mskip1mu\MA}
\def\HMAi{H\hbox{-}\mskip1mu\MAi}
\def\HMB{H\hbox{-}\mskip1mu\MB}
\def\Homk{\Hom_{\,\bbk}}
\def\HQM{H\hbox{-}\mskip2mu\vphantom{\M}_Q\M}
\def\HQMQ{\HQM_Q}
\def\Hscr{{\cal H}}
\def\IIP #1{\hbox{I\hskip1pt I\hskip1pt P}\ifx.#1\else{} \fi #1}
\def\Max{\mathop{\rm Max}\nolimits}
\def\m{\frak m}
\def\M{{\cal M}}
\def\MA{\M_A}
\def\MAi{\M_{A_i}}
\def\MB{\M_B}
\def\Nm{\mathop{\rm Nm}\nolimits}
\def\p{\frak p}
\def\q{\frak q}
\def\sHMA{\Hscr\hbox{-}\mskip1mu\sMA}
\def\sMA{\M_\Ascr}
\def\ZAcH{{Z(A_0)^{H_0}}}
\def\ZAH{{Z(A)^H}}
\def\ZQH{{Z(Q)^H}}
\def\ZQiH{{Z(Q_i)^H}}

\citedef
Art91
Art96
Bah-L98
Ber-CF90
Ber-M86
Ber-I73
Bj71
Ch90
Coh86
Coh92
Coh-F86
Coh-F92
Coh-FM90
Coh-R83
Coh-W93
Coh-W94
Coh-WZ96
De-G
Doi85
Er11
Er12
Er-a16
Er-b16
Et15
Et-W14
Et-W15
Et-W16
Fer94
Kal-T08
Kato67
Kh74
Kr-T81
Lin-M07
Matc91
McC-R
Mo78
Mo80
Mo93
Mo-S99
Mo-Sm84
Mum
Par-S78
Pas
Qu89
Row
Sch76
Sk02
Sk04
Sk07
Sk-ART11
Sk-JA11
Sk15
Sk-PAMS17
Sk-JA17
Sk-Oy06
Tot97
Tot98
Zhu96
\endcitedef

\frontmatter

\title{Subrings of invariants for actions\break
of finite dimensional Hopf algebras}
\author{Serge Skryabin}
\address{Institute of Mathematics and Mechanics,
Kazan Federal University,\break
Kremlevskaya St.~18, 420008 Kazan, Russia}
\email{Serge.Skryabin@kpfu.ru}

\endfrontmatter

The classical invariant theory on the threshold of the 20th century considered 
finite generation of invariants as a key problem. Initially only group actions 
on polynomial rings in several indeterminates were looked at. In the case of a 
finite group Emmy Noether's constructive approach via Galois resolutions 
yielded certain conclusions also in a more abstract situation. It showed that 
any commutative ring is integral over the subring of invariants with respect 
to a finite group of automorphisms. For a finitely generated commutative 
algebra over a field integrality over a subalgebra is equivalent to finiteness 
as a module over this subalgebra, while finite generation of the subalgebra is 
a consequence of these properties. This makes integrality and 
module-finiteness particularly important in the study of invariants.

Grothendieck and his school made a transition from group actions to actions of 
group schemes. As it turned out, there are general results on invariants 
incorporated in the construction of quotients by finite group schemes. These 
results can be interpreted in terms of coactions of commutative Hopf algebras 
or, dually, in terms of actions of cocommutative Hopf algebras.

In a different vein many people contributed to research on group actions and 
Lie algebra actions, and also group gradings, on noncommutative rings. The 
work on Hopf algebra actions started in the 1980s aimed to unify previously 
known results in those areas. In spite of the progress made in this study 
considerable difficulties have been encountered in some questions. Even now 
the state of knowledge in the general Hopf algebra case has not reached the 
level recorded in the 1980 manuscript of Susan Montgomery on fixed rings of 
finite automorphism groups of associative rings \cite{Mo80}.

The present paper is intended primarily as a survey of recent work on 
invariants of Hopf algebra actions. Its highlights are results on integrality 
of $H$-module PI algebras over central invariants obtained independently by 
Etingof \cite{Et15} and Eryashkin \cite{Er-a16}. There is a different 
notion of integrality introduced by Schelter \cite{Sch76} which is suitable 
for extensions of noncommutative rings. Eryashkin has also proved that an 
arbitrary $H$-module PI algebra is Schelter integral over the subring of all 
invariants when the Hopf algebra $H$ is semisimple and cosemisimple 
\cite{Er-b16}, thus answering a question of Montgomery \cite{Mo93} in the PI 
case.

Other recent results on invariants are presented in the author's two own 
papers. In \cite{Sk15} it has been proved that, given a semisimple Hopf 
algebra $H$, all nonzero $H$-stable one-sided ideals of any noetherian 
$H$-semiprime $H$-module algebra $A$ contain nonzero invariants, and the 
classical quotient ring of $A$ is obtained by localization at the Ore set of 
invariant regular elements. We will show that these conclusions are true even 
when $A$ is not noetherian provided that $A$ has an artinian classical 
quotient ring. Another article \cite{Sk-JA17} has answered the question of 
Bergen, Cohen and Fischman \cite{Ber-CF90} on the equality of the left and 
right dimensions of a skew field over the subfield of invariants. We will 
also review older results.

In this paper only finite dimensional Hopf algebras over a field will be 
considered. It should be noted, however, that many results discussed here can 
be formulated more generally when an arbitrary commutative ring is taken as a 
base ring and the Hopf algebras are finitely generated projective modules. In 
fact, this was the setting for several original papers.

\section
1. Terminology and notation

Throughout the whole paper $H$ will stand for a finite dimensional Hopf algebra 
over a field $\bbk$. We denote by $\De$, $\ep$, $S$ the comultiplication, the 
counit and the antipode in either $H$ or the dual Hopf algebra $H^*$, depending 
on the context. For a general information on Hopf algebras and their actions 
on rings we refer the reader to \cite{Art91}, \cite{Mo93} and other books.

Recall that the categories of $H$-modules and $H$-comodules are monoidal. If 
$V$ and $W$ are two (left) $H$-modules, then $V\ot W$ is an $H\ot H$-module, 
and $H$ acts on $V\ot W$ via $\De:H\to H\ot H$. If $V$ and $W$ are two (right) 
$H$-comodules, then $V\ot W$ is an $H\ot H$-comodule, and the coaction of $H$ 
is obtained by means of the map $H\ot H\to H$, $\,a\ot b\mapsto ab$. Here and 
later $\ot$ means $\ot_\bbk$ unless the base ring for the tensor product is 
indicated explicitly.

All algebras and rings are assumed to be associative and unital. An 
\emph{$H$-module algebra}\/ is a $\bbk$-algebra $A$ equipped with a left 
$H$-module structure such that the multiplication map $A\ot A\to A$ is 
$H$-linear, assuming that $H$ acts on $A\ot A$ via $\De$. If this condition is 
satisfied, then $H$ acts trivially on the image of $\bbk$ in $A$, so that 
$h1_A=\ep(h)1_A$ for all $h\in H$ where $1_A$ is the unity of $A\,$ 
\cite{Coh92, Lemma 1.9}.

The \emph{$H$-invariant}\/ elements of an $H$-module algebra $A$ form a 
subalgebra
$$
A^H=\{a\in A\mid ha=\ep(h)a\,\hbox{ for all $h\in H$}\}.
$$

An \emph{$H$-comodule algebra}\/ is a $\bbk$-algebra $A$ equipped with a right 
$H$-comodule structure such that the multiplication map $A\ot A\to A$ is 
a homomorphism of comodules. This condition can be reformulated by saying 
that the comodule structure map $\rho:A\to A\ot H$ is multiplicative, i.e. 
$$
\rho(ab)=\rho(a)\rho(b)\quad\hbox{for all $a,b\in A$}. 
$$
Moreover, $\rho$ is then a homomorphism of unital algebras. The fact that 
$\rho(1)=1\ot1$, and so $H$ coacts trivially on the image of $\bbk$ in $A$, is easily seen as follows. Clearly $\rho(1)\,x=x$ for all 
$x$ in the right ideal $I$ of $A\ot H$ generated by $\rho(A)$. So it 
suffices to check the equality $I=A\ot H$, but it does hold because the 
linear map
$$
A\ot H\to A\ot H,\qquad a\ot h\mapsto\rho(a)\cdot(1\ot h)
$$
is bijective. In fact the assignment 
$a\ot h\mapsto(\id\ot S)\bigl(\rho(a)\bigr)\cdot(1\ot h)$ defines the inverse 
map. This argument shows also that $\rho$ is an isomorphism of $A$ onto a 
subalgebra of $A\ot H$, and $A\ot H$ is a free $\rho(A)$-module with respect 
to the action by left multiplications. Similarly, $A\ot H$ is free over 
$\rho(A)$ on the right.

With an $H$-comodule algebra one associates its subalgebra consisting of 
\emph{coaction invariants}
$$
A\co H=\{a\in A\mid\rho(a)=a\ot1\}.
$$

As is well-known, the left $H$-module structures are in a bijective 
correspondence with the right $H^*$-comodule structures. This correspondence 
is compatible with the tensor products of modules and comodules. Therefore each 
$H$-module algebra is an $H^*$-comodule algebra and vice versa. Under the 
canonical identification of $A\ot H^*$ with $\Homk(H,A)$ the comodule 
structure on an $H$-module algebra $A$ is given by the map
$$
\openup1\jot
\vcenter{\halign{$#$\hfil\cr
\rho:A\to A\ot H^*\cong\Homk(H,A),\cr
\rho(a)(h)=ha\quad\hbox{for $a\in A$, $h\in H$}.\cr
}}
$$

In the rest of the paper $A$ is assumed to be an $H$-module algebra. However, 
sometimes arguments are formulated more naturally in terms of comodule 
structures. Note, in particular, that $\,A^H=A\co{H^*}$.

By an \emph{ideal}\/ we mean a two-sided ideal unless explicitly stated 
otherwise. Of particular interest are $H$-invariant ideals, i.e. ideals stable 
under the action of $H$. Several properties of an $H$-module algebra $A$ are 
defined in terms of the collection of its $H$-stable ideals:

\smallskip
$A$ is \emph{$H$-simple}\/ if $A\ne0$ and $A$ has no $H$-stable ideals except 
the zero ideal and the whole $A$;

\smallskip
$A$ is \emph{$H$-prime}\/ if $A\ne0$ and $IJ\ne0$ for all nonzero $H$-stable ideals 
$I$ and $J$ of $A$;

\medskip
$A$ is \emph{$H$-semiprime}\/ if $A$ contains no nonzero nilpotent $H$-stable 
ideals.

\medskip
An $H$-stable ideal $I$ of $A$ is called $H$-prime (respectively $H$-semiprime) 
if the factor algebra $A/I$ is $H$-prime (respectively $H$-semiprime). For an 
arbitrary ideal $I$ of $A$ we denote by $I_H$ the largest $H$-stable ideal of 
$A$ contained in $I$. If $I$ is prime (respectively semiprime), then $I_H$ is 
$H$-prime (respectively $H$-semiprime). Conversely, if $I$ is $H$-prime, then 
$I=P_H$ for some prime ideal $P$ of $A$ \cite{Ch90, Lemma 1.5}. An $H$-stable 
ideal is $H$-semiprime if and only if it is an intersection of $H$-prime 
ideals \cite{Mo-S99, Lemma 8.3}.

If a ring-theoretic notion is not prefixed by $H$-$\,$, then it does not take 
into account the $H$-module structure. For example, an $H$-module algebra $A$ 
is \emph{PI}\/ if $A$ satisfies a polynomial identity as an ordinary algebra.

A left or right $A$-module $M$ is said to be \emph{$H$-equivariant}\/ if $M$ is 
equipped with a left $H$-module structure such that the action of $A$ on $M$ 
comes from an $H$-linear map $A\ot\nobreak M\to M$ or $M\ot A\to M$, respectively, 
assuming that $H$ acts in the tensor products via $\De$.
Denote by $\HAM$ and $\HMA$ the categories of $H$-equivariant left and right 
$A$-modules. Morphisms in these categories are maps which are 
$A$-linear and $H$-linear simultaneously. Let $\AM$ and $\MA$ stand for the 
categories of left and right $A$-modules.

Similarly, an $A$-bimodule $M$ will be called $H$-equivariant if $M$ is 
equipped with a left $H$-module structure with respect to which $M$ is an 
object of both $\HAM$ and $\HMA$. We denote by $\HAMA$ the category of 
$H$-equivariant $A$-bimodules. Note that each $H$-stable ideal of $A$ is an 
object of $\HAMA$, and any homomorphism of $H$-module algebras $A\to B$ makes 
$B$ an object of $\HAMA$.

Recall that the smash product algebra $A\#H$ has $A\ot H$ as its underlying 
vector space, the canonical maps $A\to A\ot1$ and $H\to1\ot H$ are isomorphisms 
of $A$ and $H$ onto subalgebras of $A\#H$, while
$$
(1\#h)(a\#1)=\tsum h\1a\#h\2\quad\hbox{for $a\in A$, $\,h\in H$}.
$$
The compatibility of the $A$-module and $H$-module structures requested in the 
definition of the category $\HAM$ means precisely that the two module structures 
come from a single $A\#H$-module structure. Thus $\HAM$ is identified with the 
category of left $A\#H$-modules.

Let $A\op$ be $A$ with the opposite multiplication, and let $H\cop$ be $H$ 
with the opposite comultiplication. Then $A\op$ is an $H\cop$-module algebra, 
and $\HMA$ may be identified with the category of left $A\op\#H\cop$-modules.

The algebras $A\#H$ and $A^H$ are connected by a Morita context. Several 
articles \cite{Ber-CF90}, \cite{Ber-M86}, \cite{Coh-F86}, \cite{Coh-FM90} 
derive various information about the invariant ring $A^H$ when $A\#H$ is known 
to be simple, or prime, or semiprime. These results are also discussed in 
\cite{Mo93}. However, for an arbitrary finite dimensional Hopf algebra $H$ it 
is quite difficult to understand the ring structure of $A\#H$ in terms of the 
original algebra $A$. There is still a big gap between what is known in the 
general case and in the case of a finite group $G$ acting on a ring $R$ where 
the skew group ring $R*G$ is sufficiently well understood as a normalizing 
extension of $R$. In our paper we rarely use ring-theoretic properties of 
$A\#H$ directly. However, equivariant modules are important for many 
considerations.

Recall that a \emph{left}\/ (respectively \emph{right}\/) \emph{integral}\/ in 
$H$ is an element $0\ne t\in H$ such that $ht=\ep(h)t$ (respectively 
$th=\ep(h)t$\/) for all $h\in H$. Let $t$ be a left integral. If $V$ is a left 
$H$-module, then the action of $t$ gives a map $\hatt:V\to V^H$ where $V^H$ is 
the subspace of $H$-invariant elements in $V$. In \cite{Coh-FM90} this map was 
called a \emph{trace}\/ by analogy with the terminology used in the case of 
group actions.

By Maschke's theorem $H$ is semisimple if and only if $\ep(t)\ne0$. In this 
case $t$ acts on $V^H$ as a nonzero scalar multiplication. In particular, 
we have $tV=V^H$, i.e. the trace $\hatt$\/ is always surjective.

\section
2. Structural properties of $H$-module algebras

In this section we present several results concerned with the structure of 
\hbox{$H$-module} algebras on which recent work on invariants of Hopf actions 
is based. Actually most of these results can be formulated for a not necessarily 
finite dimensional Hopf algebra $H$. Nevertheless it will be assumed in all 
statements that $\,\dim H<\infty$. With this assumption we do not need to 
mention any additional restrictions, and also the proofs become considerably 
simpler.

A key argument used in deriving these results comes from

\proclaim
Theorem 2.1 \cite{Sk07}.
Suppose that $A$ is a semilocal $H$-simple $H$-module algebra. Then each 
object $M\in\HMA$ is projective in $\MA$. Moreover, a direct sum of several 
copies of $M$ is a free $A$-module. A similar conclusion holds in $\HAM$.
\endproclaim

There is one application of the freeness properties of $H$-equivariant modules 
where the $H$-simplicity of the $H$-module algebra is not known beforehand. To 
deal with this situation one needs the next lemma stated under more technical 
assumptions about $A$ and $M$ than the previous theorem.

Denote by $\Max A$ the set of maximal ideals of $A$. If $A$ is semilocal, then 
its factor algebra by the Jacobson radical is semisimple artinian. This means 
that the set $\Max A$ is finite and $A/P$ is simple artinian for each 
$P\in\Max A$. An object $M$ of the category $\HMA$ is said to be 
\emph{$A$-finite}\/ if $M$ is finitely generated as an $A$-module. If $M$ is 
$A$-finite, then $M/MP$ is an $A$-module of finite length. The \emph{rank}\/ 
of $M$ at $P$ is defined 
as
$$
r_P(M)={\lng_AM/MP\over\lng_AA/P}\in\bbQ\,.
$$

\proclaim
Lemma 2.2.
Let $A$ be a semilocal $H$-module algebra. Suppose that $M\in\HMA$ is $A$-finite 
and there exists $P\in\Max A$ such that $P$ contains no nonzero $H$-stable 
ideals of $A$ and $r_P(M)\ge r_Q(M)$ for all $Q\in\Max A$. Then $M^n$ is a free 
$A$-module for some integer $n>0$.
\endproclaim

For the proof see \cite{Sk07, Lemma 7.5}. This lemma is valid even when $H$ is 
an infinite dimensional Hopf algebra. However, the assumption $\dim H<\infty$ 
is needed to deduce the conclusion of Theorem 2.1 for objects $M$ which are not 
$A$-finite. Since each element of $M$ is contained in an $A$-finite subobject 
of $M$, a basis for $M^n$ over $A$ can be constructed using Zorn's lemma.

It turns out that Lemma 2.2 and Theorem 2.1 lead to several fundamental facts 
concerning artinian $H$-module algebras very quickly. Sometimes one has 
initially less information about an $H$-module algebra, but the left and right 
artinian conditions can be deduced. For this reason we have to deal with 
semiprimary algebras. A semilocal ring is called \emph{semiprimary}\/ if its 
Jacobson radical is nilpotent.

\proclaim
Lemma 2.3.
Let $A$ be a semiprimary $H$-module algebra, and let $K\in\Max A$. Then the 
largest $H$-stable ideal $K_H$ contained in $K$ is a maximal $H$-stable 
ideal of $A$.
\endproclaim

\Proof.
Replacing $A$ with the factor algebra $A/K_H$, we may assume that $K_H=0$, and 
then we have to prove that $A$ is $H$-simple. First note that $A$ is $H$-prime. 
Indeed, if $I,J$ are two nonzero $H$-stable ideals of $A$, then both 
$I\not\sbs K$ and $J\not\sbs K$, whence $IJ\not\sbs K$, and therefore $IJ\ne0$.

Every semiprimary ring satisfies DCC on finitely generated one-sided ideals. 
Hence $A$ has a minimal nonzero $H$-stable finitely generated right ideal $M$.  
If $0\ne x\in M$, then $M=(Hx)A$ since $(Hx)A$ is a nonzero $H$-stable 
finitely generated right ideal of $A$ contained in $M$. It follows that $M$ is 
minimal in the set of all nonzero $H$-stable right ideals of $A$. If $I$ is 
any nonzero $H$-stable ideal of $A$, then $MI$ is an $H$-stable right ideal. 
Since $MI\sbs M$ and $MI\ne0$ by the $H$-primeness of $A$, we get $MI=M$.

We may view $M$ as an $A$-finite object of $\HMA$. Pick $P\in\Max A$ for which 
$r_P(M)$ attains the maximum value. Since $M\ne0$, we have $r_P(M)>0$. This 
means that $M\ne MP$, but then $M\ne MP_H$ too, which is only possible when 
$P_H=0$ by the preceding argument. Thus the hypotheses of Lemma 2.2 are 
fulfilled, and we deduce that $M^n$ is a free $A$-module for some $n>0$. Hence 
$MI\ne M$ for each ideal $I\ne A$. If $I$ is $H$-stable and $I\ne A$, this 
entails $I=0$.
\endproof

\proclaim
Theorem 2.4 \cite{Sk-Oy06}.
Suppose that $A$ is semiprimary and $H$-semiprime. There is an isomorphism of 
$H$-module algebras
$$
A\cong A_1\times\ldots\times A_n
$$ where $A_1,\ldots,A_n$ 
are $H$-simple $H$-module algebras.
If $A$ has a maximal ideal containing no nonzero $H$-stable ideals of $A$ then 
$A$ is $H$-simple.
\endproclaim

\Proof.
By Lemma 2.3 the maximal $H$-stable ideals of $A$ are precisely the ideals $K_H$ 
with $K\in\Max A$. In particular, there are finitely many of them. Let 
$I_1,\ldots,I_n$ be all the maximal $H$-stable ideals. Then $I_1\cap\ldots\cap 
I_n$ is contained in the Jacobson radical $J$ of $A$. Since $J$ is nilpotent 
and $A$ is $H$-semiprime, we get $I_1\cap\ldots\cap I_n=0$. But $I_k+I_l=A$ 
for each pair of indices $k\ne l$, whence the desired direct product 
decomposition of $A$ holds with $A_k=A/I_k$ by the Chinese remainder theorem.  
\endproof

\proclaim
Corollary 2.5.
Each right coideal subalgebra $B$ of $H^*$ is an $H$-simple $H$-module algebra 
and $H^*$ is right and left $B$-free.
\endproclaim

\Proof.
Here $B$ is a subalgebra and a right coideal of $H^*$. The restriction of the 
comultiplication $\De$ in $H^*$ gives a map $B\to B\ot H^*$ which makes $B$ 
into an $H^*$-comodule algebra. Hence $B$ is also an $H$-module algebra. Put 
$B^+=\Ker\ep|_B$ where $\ep$ is the counit of $H^*$. Then $B^+\in\Max B$ with 
$B/B^+\cong\bbk$. Suppose that $I$ is an $H$-stable ideal of $B$. Then 
$\De(I)\sbs I\ot H^*$. Recall that $(\ep\ot\id)\circ\De$ is the identity map 
by the definition of the counit. If $I\sbs B^+$, we get $x=(\ep\ot\id)(\De 
x)=0$ for each $x\in I$ since $\ep(I)=0$, and so $I=0$. Thus $B^+$ contains no 
nonzero $H$-stable ideals of $B$. By Theorem 2.4 $B$ is $H$-simple.

Now $H^*$ is also an $H$-module algebra, and $B$ is its $H$-stable subalgebra.  
Hence we may regard $H^*$ as an object of both $\HBM$ and $\HMB$. Freeness of 
$H^*$ over $B$ follows from Theorem 2.1.
\endproof

\proclaim
Corollary 2.6.
Suppose that $A$ is semiprimary and $H$-semiprime. Then each object $M\in\HMA$ 
is projective in $\MA$.
\endproclaim

\Proof.
The direct product decomposition of $A$ given in Theorem 2.4 implies that 
$M\cong M_1\times\ldots\times M_n$ where $M_i=M\ot_AA_i\in\HMAi$. By Theorem 
2.1 $M_i$ is projective in $\MAi$ for each $i$, whence the conclusion.
\endproof

\proclaim
Corollary 2.7.
Suppose that $A$ is semiprimary and $H$-semiprime. If $I$ is an $H$-stable 
right ideal of $A$, then $I=eA$ for some idempotent $e\in A$.
\endproclaim

\Proof.
We may view $A/I$ as an object of $\HMA$. By Corollary 2.6 $A/I$ is projective 
in $\MA$. Hence $I$ is a direct summand of $A$ as a right $A$-module.
\endproof

\proclaim
Theorem 2.8 \cite{Sk-Oy06}.
Any semiprimary $H$-semiprime algebra $A$ is a quasi-Frobenius ring. In 
particular, $A$ is left and right artinian.
\endproclaim

\Proof.
By Theorem 2.4 it suffices to consider the case when $A$ is $H$-simple. We are 
going to apply the general fact that a semiprimary ring is quasi-Frobenius 
whenever it is left and right selfinjective (see \cite{Kato67, Th. 10}). Let 
us show that $A$ is left selfinjective. Applying this to the $H\cop$-module 
algebra $A\op$, we deduce that $A$ is right selfinjective too, and the 
conclusion follows.

Take any nonzero injective object $M\in\HAM$. Then $M$ remains injective in 
$\AM$. To see this recall that $\HAM$ is identified with the category of left 
$B$-modules for $B=A\#H$. The forgetful functor $\HAM\to\AM$ is identified 
with the restriction functor $\BM\to\AM$ that arises from the canonical 
embedding of $A$ into $B$. The latter functor preserves injectives since it 
has an exact left adjoint $\,B\,\ot_A\,?\,$.

But $M^n$ is a free $A$-module for some $n>0$ by Theorem 2.1. Therefore $A$ is an 
$\AM$-direct summand of $M^n$. Since $M^n$ is injective in $\AM$, so too is $A$.
\endproof

Thus an $H$-semiprime algebra $A$ is semiprimary if and only if $A$ is left 
artinian, if and only if $A$ is right artinian, and we will say that $A$ is 
\emph{artinian}\/ in this case.

\proclaim
Theorem 2.9.
Let $A$ be an $H$-stable subalgebra of an $H$-module algebra $B$. If $A$ is 
artinian and $H$-simple then $B$ is free as an $A$-module with respect to the 
action by right {\rm(}or left{\rm)} multiplications.
\endproclaim

\Proof.
We may view $B$ as an object of $\HMA$. Hence Theorem 2.1 applies to $B$. For 
each $P\in\Max A$ denote by $F_P$ the projective cover in $\M_A$ of a simple 
right $A/P$-module. These modules $F_P$ are indecomposable, and 
$A\cong\bigoplus_{P\in\Max A}F_P^{\,m_P}$ for some multiplicities $m_P$. Put
$$
\hbox{$E=\bigoplus_{P\in\Max A}F_P^{\,m_P/d}$\quad
where $d=\gcd\{m_P\mid P\in\Max A\}$}.
$$
Then $A\cong E^d$. If $M$ is any right $A$-module such that $M^n$ is free for 
some $n>0$, then, by the Krull-Schmidt theorem, $M$ is isomorphic to a direct 
sum of a family of copies of $E$. Moreover, $M$ is itself free if either $M$ 
is not finitely generated or $M\cong E^k$ with $d$ dividing $k$.
If $M$ is in fact an $A$-bimodule, then $M\cong A\ot_AM\cong N^d$ where 
$N=E\ot_AM$. In this case $N$ has to be isomorphic to a direct sum of a family 
of copies of $E$, whence $M$ is right $A$-free by the previous observation. It 
remains to apply this for $M=B$.
\endproof

\proclaim
Theorem 2.10 \cite{Sk-Oy06}.
Suppose that $H$ is semisimple, $A$ is artinian and $H$-semiprime. Then 
$\HMA$ and $\HAM$ are semisimple categories. In other words, the smash 
product algebras $A\op\#H\cop$ and $A\#H$ are semisimple artinian.
\endproclaim

\Proof.
A key ingredient in the proof is the fact established by Cohen and Fischman 
\cite{Coh-F86} according to which a submodule $W$ of a left $A\#H$-module $V$ 
is a direct summand whenever $W$ is an $A$-module direct summand of $V$. For 
this one needs only semisimplicity of $H$ but no assumptions about an 
$H$-module algebra $A$.

In the case of the category $\HMA$ a similar argument runs as follows. Let 
$M,N$ be two objects of $\HMA$. There is a left $H$-module structure on 
$\Homk(M,N)$ defined by the rule
$$
(h\cdot f)(x)=\sum h\1\,f\bigl(S(h\2)x\bigr),\qquad
$$
for $h\in H$ with $\De h=\sum h\1\ot h\2$, $\,f\in\Homk(M,N)$, $\,x\in M$.
It is straightforward to check that $\Hom_A(M,N)$ is stable under this action 
of $H$ and that a linear map $f\in\Homk(M,N)$ is $H$-invariant if and only 
if $f$ is $H$-linear. Let $t\in H$ be an integral with $\ep(t)=1$. If $N$ is 
a subobject of $M$ which splits off as an $A$-module direct summand, then 
there exists an $A$-linear map $f:M\to N$ such that $f|_N=\id$. Now the map 
$f'=t\cdot f$ is $A$-linear and $H$-linear simultaneously, and also 
$f'|_N=\id$. Hence $M=N\oplus\Ker f'$, a direct sum decomposition in $\HMA$.

Under the assumption that $A$ is artinian and $H$-semiprime any subobject $N$ 
of $M$ is an $A$-module direct summand since the factor object $M/N\in\HMA$ is 
projective in $\MA$ by Corollary 2.6. Hence the previous conclusion holds.
\endproof

The result of Cohen and Fischman mentioned in the proof of Theorem 2.10 means 
that $A\#H$ is a semisimple extension of $A$ when $H$ is semisimple. In 
\cite{Coh-F86} it was used to show that $A\#H$ is semiprime artinian whenever 
so is $A$. If $A$ is not semiprime but only $H$-semiprime, the same conclusion 
requires Theorem 2.1 which has been proved much later.

Let $R$ be a ring. A ring $Q$ is said to be a \emph{classical right quotient 
ring}\/ of $R$\/ if

\setitemsize(3)
\item(1)
$Q$ contains $R$ as a subring,

\item(2)
all \emph{regular}\/ elements, i.e. nonzerodivisors, of $R$ are invertible in 
$Q$, and

\item(3)
each element $q\in Q$ can be written as $q=as^{-1}$ where $a,s\in R$, $\,s$ 
regular.

\noindent
Such a ring $Q$ exists if and only if the set of all regular elements of $R$ 
satisfies the right Ore condition. In this case $Q$ is unique up to 
isomorphism, and we will denote this ring by $Q(R)$.

If an $H$-module algebra $A$ is not artinian, but $A$ has an artinian classical 
right quotient ring, the previous results can still be used to derive 
information about $A$. There are two important cases when this happens:

\proclaim
Theorem 2.11 \cite{Sk-Oy06}.
If $A$ is right noetherian and $H$-semiprime, then $A$ has a quasi-Frobenius 
classical right quotient ring.
\endproclaim

\proclaim
Theorem 2.12 \cite{Er-a16}.
If $A$ is PI and $H$-semiprime with finitely many minimal $H$-prime ideals, 
then $A$ has a quasi-Frobenius classical right quotient ring. In particular, 
this holds if $A$ is finitely generated, PI, and $H$-semiprime.
\endproclaim

In the proof of Theorem 2.11 one first constructs a generalized quotient ring 
$Q$ using the filter of $H$-stable essential right ideals of $A$. This ring 
turns out to be semiprimary. In the proof of Theorem 2.12 one starts with the 
$H$-equivariant Martindale quotient ring $Q$. As we shall see in section 6 it 
will be a finite module over a central artinian subring. In both cases $H$ 
acts on $Q$, and $H$-semiprimeness is preserved under passage to $Q$. Hence 
$Q$ is quasi-Frobenius by Theorem 2.8, and the conclusion that $Q$ is a 
classical right quotient ring can be deduced from the following ring-theoretic 
fact:

\proclaim
Proposition 2.13 \cite{Sk-JA11}.
Let $R$ be a subring of a quasi-Frobenius ring $Q$. Suppose that $\I$ is a 
topologizing filter of right ideals of $R$ with the following properties\/{\rm:}

\item(a)
each $I\in\I$ has zero left and right annihilators in $Q${\rm,}

\item(b)
for each $q\in Q$ there exists $I\in\I$ such that $qI\sbs R$.

\noindent
Then each right ideal $I\in\I$ contains a regular element of $R$, and $Q$ is a 
classical right quotient ring of $R$.
\endproclaim

We say that a family $\F$ of right ideals in a ring $R$ is a {\it filter}\/ 
if for any pair of right ideals $I,J\in\F$ there exists $K\in\F$ such that 
$K\sbs I\cap J$. A filter $\F$ is \emph{topologizing}\/ if for each $I\in\F$ 
and each $a\in R$ there exists $I'\in\F$ such that $aI'\sbs I$. The condition 
that with each $I\in\F$ all larger right ideals also belong to $\F$ is often 
included in the definition of a filter, but omitting it will do no harm.

With small improvements in the proof of \cite{Sk-JA11, Prop. 1.4} the assumption 
that the filter $\I$ is topologizing can be actually removed. In the case when 
$Q$ is semisimple artinian see \cite{Sk15, Prop. 2.3}.

\proclaim
Theorem 2.14 \cite{Sk-Oy06}.
Suppose that $A$ has a right artinian classical right quotient ring $Q$. Then
the $H$-module structure on $A$ has a unique extension to $Q$ with respect to
which $Q$ becomes an $H$-module algebra.
\endproclaim

\Proof.
We argue in terms of comodule structures. Since the map
$$
A\ot H^*\to A\ot H^*,\qquad a\ot f\mapsto\rho(a)\cdot(1\ot f),
$$
is invertible, $A\ot H^*$ is a free $A$-module with respect to the action of 
$A$ given by left multiplications by the elements $\rho(a)$. For each regular 
element $s$ of $A$ it follows that $\rho(s)$ is right regular in $A\ot H^*$, 
i.e. $\rho(s)x=0$ for $x\in A\ot H^*$ implies $x=0$. Then $\rho(s)$ remains 
right regular in $Q\ot H^*$, and therefore $\rho(s)$ has to be invertible in 
the right artinian ring $Q\ot H^*$. This property shows that $\rho:A\to A\ot 
H^*$ extends to an algebra homomorphism $\rho':Q\to Q\ot H^*$. Now 
$(\id\ot\De)\rho'$ and $(\rho'\ot\id)\rho'$ are two algebra homomorphisms 
$Q\to Q\ot H^*\ot H^*$ which agree on $A$. Hence
$$
(\id\ot\De)\rho'=(\rho'\ot\id)\rho'.
$$
Thus $\rho'$ is a structure of an $H^*$-comodule algebra extending the given 
one on $A$.
\endproof

In the situation of Theorem 2.14 we have two subrings of invariants $A^H$ 
and $Q^H$. Clearly $A^H=A\cap Q^H$. Theorem 2.14 is true even without the 
assumption that $\dim H<\infty$.

The conclusions of Theorems 2.4, 2.8 and 2.11 hold for some classes of 
infinite dimensional Hopf algebras too. Unfortunately it is still unknown 
whether the assumption that $H$ has a bijective antipode is sufficient for 
their validity. However, if $A$ is right artinian and $H$-semiprime, then $A$ 
is a quasi-Frobenius ring, even when $H$ is an arbitrary infinite dimensional 
Hopf algebra \cite{Sk-ART11}.

\section
3. Module-finiteness over the invariants

In this section we examine those cases where $A$ is known to be a finite 
$A^H$-module. Many results discussed here date back to as early as the 1980s.

The comodule structure $\rho:A\to A\ot H^*$ enables us to define a $\bbk$-linear 
map
$$
\ga:\AHA\to A\ot H^*,\qquad
a\ot b\mapsto(a\ot1)\cdot\rho(b)
$$
(recall that $A^H=\{a\in A\mid\rho(a)=a\ot1\}$). Under the canonical 
identification of $A\ot H^*$ with $\Homk(H,A)$ we have
$$
\ga(a\ot b)(h)=a(hb)\quad\hbox{ for $a,b\in A$ and $h\in H$}.
$$
The notion of Hopf Galois extensions of algebras is defined in terms of comodule 
structures (see \cite{Mo93, Ch. 8}). In the case when the Hopf algebra is 
finitely generated projective as a module over a commutative base ring this 
notion was introduced by Kreimer and Takeuchi \cite{Kr-T81}.

Since $A$ is an $H$-module algebra by our convention, we say that $A$ is an 
\emph{$H^*$-Galois extension} of the subalgebra $A^H$ if $\ga$ is bijective. 
Since $\dim H^*<\infty$, it suffices to require surjectivity of $\ga\,$ 
\cite{Kr-T81}.

\proclaim
Theorem 3.1 \cite{Kr-T81}.
Suppose that $A$ is an $H^*$-Galois extension of $A^H$. Then $A$ is a finitely 
generated projective $A^H$-module on the left and on the right.
\endproclaim

\Proof.
As explained in \cite{Mo93, 8.3.1} this conclusion can be proved by verifying 
the dual basis property which characterizes projective modules. To do this 
let $t\in H$ be a left integral and $\la\in H^*$ a right integral such that 
$\la(t)=1$, and let $1\ot\la=\ga(\,\sum a_i\ot b_i)$ for some elements 
$a_i,b_i\in A\,$ ($i=1,\ldots,n$). If $x\in A$, then
$$
\ga(\,\tsum a_i\ot b_ix)=(1\ot\la)\cdot\rho(x)=x\ot\la,
$$
$$
\hbox{which means that}\quad
\tsum a_i\bigl(h(b_ix)\bigr)=\la(h)x\quad
\hbox{for all $h\in H$}.
$$
For each $i$ define a right $A^H$-linear map $f_i:A\to A^H$ by the formula 
$f_i(x)=t(b_ix)$. Taking $h=t$ above, we deduce that
$$
\tsum a_if_i(x)=\tsum a_i\bigl(t(b_ix)\bigr)=\la(t)x=x
$$
for all $x\in A$. In the proof of the left side version one proceeds similarly 
replacing $\ga$ with the map $\,a\ot b\mapsto\rho(a)\cdot(b\ot1)\,$ which is 
also bijective by \cite{Kr-T81, Prop. 1.2}.
\endproof

Hopf Galois extensions form an important special class of comodule algebras 
for which more information is available than in general. However, the map 
$\ga$ is quite useful, even when $\ga$ is not bijective. Let us view 
$A\ot H^*$ as an $A$-bimodule with respect to the left and right actions 
defined by the rules 
$$
ax=(a\ot1)\cdot x,\qquad xa=x\cdot\rho(a)\quad
\hbox{for $a\in A$ and $x\in A\ot H^*$}.
$$
Then $\ga$ is a homomorphism of $A$-bimodules. In particular, its image is a 
subbimodule of $A\ot H^*$. Also, $\ga$ respects certain $H$-module structures. 
Recall the two natural left actions of $H$ on $H^*$ defined by the formulas
$$ 
(h\rhu\xi)(g)=\xi(gh),\qquad(h\rhd\xi)(g)=\xi\bigl(S(h)g\bigr)\quad
\hbox{for $g,h\in H$ and $\xi\in H^*$}.
$$

\proclaim
Lemma 3.2.
With respect to the $H$-module structures on $\AHA$ and $A\ot H^*$ defined by 
the formulas
$$
h(a\ot b)=ha\ot b,\qquad h(a\ot\xi)=\tsum h\1a\ot(h\2\rhd\xi)\eqno(*)
$$
where $a,b\in A$, $\,h\in H$ and $\xi\in H^*$ the map $\ga$ is a morphism in 
$\HAM$. On the other hand, $\ga$ is a morphism in $\HMA$ with respect to 
another pair of $H$-module structures
$$
h(a\ot b)=a\ot hb,\qquad h(a\ot\xi)=a\ot(h\rhu\xi).\eqno(**)
$$
\endproclaim

\Proof.
Clearly the $H$-module structures $(*)$ are compatible with the left 
$A$-module structures, so that $\AHA$ and $A\ot H^*$ become objects of $\HAM$. 
To show that $\ga$ is $H$-linear with respect to $(*)$ we need only to check 
that $\rho(b)$ is $H$-invariant for each $b\in A$. In the monoidal category of 
left $H$-modules the module $H^*$ with the action $\rhd$ is the left dual of 
$H$ with the action by left multiplications. Therefore
$$
\Hom_H(V,A\ot H^*)\cong\Hom_H(V\ot H,A)
$$
for each left $H$-module $V$. Taking $V=\bbk$ with the trivial module structure, 
we get $(A\ot H^*)^H\cong\Hom_H(H,A)$. Under this bijection the $H$-linear map 
$H\to A$, $\,h\mapsto hb$, corresponds to $\rho(b)\in A\ot H^*$. In other words, 
$(A\ot H^*)^H=\rho(A)$.

Now $H^*$ is an $H$-module algebra with respect to the action $\rhu$. Hence so 
too is $A\ot H^*$ with respect to the action given in $(**)$. The map 
$\rho:A\to A\ot H^*$ is $H$-linear with respect to this action, and so $\rho$ 
is a homomorphism of $H$-module algebras. In particular, $A\ot H^*\in\HMA$. 
Clearly $\AHA\in\HMA$ too. Finally, the map $\ga$ is $H$-linear with respect 
to $(**)$ since the elements in $A\ot1$ are $H$-invariant.
\endproof

\proclaim
Lemma 3.3.
Put $M=\Im\ga$. Suppose that $M$ is left $A$-free and there exist elements 
$a_1,\ldots,a_n\in A$ such that $\rho(a_1),\ldots,\rho(a_n)$ form a basis for 
$M$ over $A\ot1$. Then $a_1,\ldots,a_n$ are a basis for $A$ as an $A^H$-module 
with respect to the action by left multiplications. In particular, $A$ is left 
free of finite rank over $A^H$.
\endproclaim

\Proof.
Given $a\in A$, there exist uniquely determined elements $c_1,\ldots,c_n\in A$ 
such that
$$
\rho(a)=\tsum\,(c_i\ot1)\rho(a_i)=\ga(\,\tsum c_i\ot a_i). 
$$
Then $ha=\tsum c_i(ha_i)$ for all $h\in H$. Taking $h=1$, we get $a=\sum c_ia_i$. 
Now let $g\in H$. Since $\rho(a)$ is $H$-invariant and $\ga$ is $H$-linear with 
respect to the $H$-module structures $(*)$ considered in Lemma 3.2, we deduce that
$$
\tsum\,(gc_i\ot1)\rho(a_i)=\ga(\,\tsum gc_i\ot a_i)=\ep(g)\rho(a)
=\ep(g)\tsum\,(c_i\ot1)\rho(a_i).
$$
It follows that $gc_i=\ep(g)c_i$ since $\rho(a_1),\ldots,\rho(a_n)$ are left 
linearly independent over $A\ot1$. Hence $c_i\in A^H$ for each $i$. 
On the other hand, if $a=\sum c'_ia_i$ for another collection of elements 
$c'_1,\ldots,c'_n\in A^H$, then $\rho(a)=\sum\,(c'_i\ot1)\rho(a_i)$, 
which entails $c'_i=c_i$ for each $i$.
\endproof

If $A$ is artinian and $H$-simple, then the $A$-bimodule $M=\Im\ga$ is always 
left (and right) free. Indeed, $M$ may be regarded as an object of $\HAM$ by 
Lemma 3.2. Hence $M^n$ is left free for some $n>0$ by Theorem 2.1. Moreover, 
this conclusion holds with $n=1$, as explained in the proof of Theorem 2.9. As 
a left $A$-module, $M$ is generated by $\rho(A)$, but this does not mean that a 
basis can be chosen in $\rho(A)$, and so Lemma 3.3 does not apply in general. 
Here lies the source of possible misbehavior of the subring $A^H$.

An example of Bj\"ork \cite{Bj71} produces a simple artinian ring $R$ of 
characteristic 2 which is not a finitely generated module over the subring 
$R^G$ of elements fixed by an automorphism group of order 4. In this example 
$R^G$ is neither left nor right artinian.

The situation becomes much nicer when the $H$-module algebra $A$ has no 
nontrivial $H$-stable one-sided ideals.

\proclaim
Lemma 3.4.
Suppose that $A$ has no nontrivial $H$-stable left {\rm(}or right{\rm)} ideals. 
Then $A^H$ is a skew field and $\ga$ is injective.
\endproclaim

\Proof.
The condition imposed on $A$ means that $A$ is a simple object of $\HAM$, i.e. 
a simple left $A\#H$-module. The fact that $A^H$ is a skew field is stated 
in \cite{Ber-CF90, Lemma 2.1}. It follows from Schur's lemma since 
$\,A^H\cong(\End_{A\#H}A)\op$.

As explained in Lemma 3.2, $\AHA$ may be regarded as an object of $\HAM$. It 
is a sum of its simple subobjects $A\ot b$ with $0\ne b\in A$, each isomorphic 
to $A$. Hence any subobject of $\AHA$ in the category $\HAM$ is equal to 
$A\ot_{A^H}\!V$ for some left vector subspace $V$ of $A$ over $A^H$. In 
particular, this applies to the kernel of $\ga$. But $\ga(1\ot b)=\rho(b)$ 
for each $b\in A$. This shows that the restriction of $\ga$ to $1\ot A$ is 
injective, and therefore $\,\Ker\ga=0$.
\endproof

Under the hypothesis of Lemma 3.4 the $H$-module algebra $A$ may be regarded 
as either left or right vector space over the skew field $A^H$. Denote by 
$[A:A^H]_l$ and $[A:A^H]_r$ the dimensions of these two vector spaces.

\proclaim
Theorem 3.5 \cite{Ber-CF90}.
Suppose that $A$ has no nontrivial $H$-stable left ideals and that $A$ has 
finite left Goldie rank {\rm(}i.e., $A$ satisfies ACC on direct sums of left 
ideals{\rm)}. Then $[A:A^H]_r\le n$ where $n$ is the dimension of the image of 
$H$ in\/ $\Endk A$.
\endproclaim

This result of Bergen, Cohen and Fischman \cite{Ber-CF90, Th. 2.2} was stated 
for a not necessarily finite dimensional Hopf algebra $H$ with finite 
dimensional image $\pi(H)$ in $\Endk A$. In the proof given in \cite{Ber-CF90} 
an application of Jacobson's density theorem shows that, whenever $A$ contains 
$m$ right linearly independent over $A^H$ elements, there exists a free left 
$A$-submodule of rank $m$ in the image $\pi(A\#H)$ of $A\#H$ in $\Endk A$. 
Pick $h_1,\ldots,h_n\in H$ whose images give a basis for $\pi(H)$. Then 
$\pi(A\#H)=\pi(T)$ where $T\sbs A\#H$ is the left $A$-submodule generated by 
$1\#h_1,\ldots,1\#h_n$. It follows that $T$ contains a free $A$-submodule 
of rank $m$. Since $T$ is a free $A$-module of rank $n$, this entails $m\le n$ 
by finiteness of the Goldie rank.

Using the map $\ga$ we can strengthen the previous theorem. We continue to work 
under the assumption that $\,\dim H<\infty$. Note, however, that replacing 
$A\ot H^*$ with $\Homk(H,A)$ in the preceding discussion and modifying all 
arguments appropriately, the next result can be proved for an infinite 
dimensional Hopf algebra under the assumption $\,\dim\pi(H)<\infty\,$ used in 
\cite{Ber-CF90}. In the semilocal case one needs also bijectivity of the 
antipode.

\proclaim
Theorem 3.6.
Suppose that $A$ has no nontrivial $H$-stable left ideals and that either 
$A$ has finite left Goldie rank or $A$ is semilocal. Let $n=\dim\pi(H)$ where 
$\pi:H\to\Endk A$ is the representation afforded by the action of $H$ on $A$. 
Then
$$
[A:A^H]_l\le n\quad\hbox{and\/}\quad[A:A^H]_r\le n.
$$
In particular, $A$ is left and right artinian.
\endproclaim

\Proof.
By Lemma 3.4 $\ga$ is injective. Thus the $A$-bimodule $M=\Im\ga$ is isomorphic 
to $\AHA$. Hence $M$ is left free of rank equal to $[A:A^H]_l$, and $M$ is right 
free of rank equal to $[A:A^H]_r$. Note that $\rho(A)\sbs A\ot C$ where $C$ is 
the subcoalgebra of $H^*$ dual to the factor algebra $\pi(H)$ of $H$. Hence 
$M\sbs A\ot C$ too. Since $\dim C=n$, this shows that $M$, regarded as an 
$A$-module with respect to the left action, embeds in a free $A$-module of rank 
$n$. If $A$ has finite left Goldie rank, we get $[A:A^H]_l\le n$, while the 
second inequality is the content of Theorem 3.5.

Suppose further that $A$ is semilocal. Since $C$ is stable under the action 
$\rhd$ of $H$ on $H^*$, the left $A$-module $A\ot C$ is an object of $\HAM$ 
with respect to the $H$-module structure described in $(*)$ of Lemma 3.2. 
Hence $(A\ot C)/M\in\HAM$ too. Since $A$ is $H$-simple, Theorem 2.1 shows that 
$(A\ot C)/M$ is projective in $\AM$. Then $M$ is an $\AM$-direct summand of 
$A\ot C$. Denoting by $J$ the Jacobson radical of $A$, we deduce that the free 
$A/J$-module $M/JM$ of rank equal to $[A:A^H]_l$ embeds into the free 
$A/J$-module $A/J\ot C$ of rank $n$. Since $A/J$ is semisimple artinian, we 
must have $[A:A^H]_l\le n$.

For the remaining part put $N=\bigl(1\ot S^{-1}(C)\bigr)\cdot\rho(A)$. Then 
$N$ is a subobject of $A\ot H^*$ in the category $\HMA$, and $N$ is $A$-free 
of rank $n$. Let $a\in A$. Writing symbolically $\rho(a)=\sum a\0\ot a\1\in 
A\ot C$, we have
$$
a\ot 1=\tsum a\0\ot S^{-1}(a\2)a\1
=\tsum\,\bigl(1\ot S^{-1}(a\1)\bigr)\cdot\rho(a\0)\in N.
$$
This shows that $A\ot1\sbs N$. Hence $M=(A\ot1)\,\rho(A)\sbs N$. Applying 
Theorem 2.1 to the object $N/M\in\HMA$, we deduce that $M$ is an $\MA$-direct 
summand of $N$, and passing to quotients modulo $J$, we arrive at 
$\,[A:A^H]_r\le n$.
\endproof

Another old result of Cohen, Fischman and Montgomery determines when the 
extension $A/A^H$ is $H^*$-Galois for an $H$-module algebra satisfying the 
previous assumptions. We will show how the map $\ga$ can be used in the proof.

\proclaim
Theorem 3.7 \cite{Coh-FM90}.
Suppose that $A$ has no nontrivial $H$-stable left ideals and that $A$ has 
finite left Goldie rank. Then the following conditions are equivalent\/{\rm:}

\setitemsize(1')
\item(1)
$[A:A^H]_r=\dim H$.

\item($1'$)
$[A:A^H]_l=\dim H$.

\item(2)
$A$ is a faithful left $A\#H$-module.

\item(3)
$A\#H$ is a simple algebra.

\item(4)
$A$ is an $H^*$-Galois extension of $A^H$.

\endproclaim

\Proof.
The map $\ga$ is an injective homomorphism of $A$-bimodules. Both $\AHA$ and 
$A\ot H^*$ are free of finite rank as left $A$-modules and as right ones. 
Since $A$ is left and right artinian by Theorem 3.6, all finitely generated 
$A$-modules have finite length. Therefore $\ga$ is bijective if and only if 
the two bimodules have equal left ranks, and if and only if they have equal 
right ranks. This shows that $(4)\Lrar(1')\Lrar(1)$.

Denote by $E$ the endomorphism ring of $A$ as a right vector space over the 
skew field $A^H$. The action of $A\#H$ on $A$ gives a ring homomorphism 
$\pi:A\#H\to E$. Since $[A:A^H]_r<\infty$, the ring $E$ is simple artinian. 
By one of general characterizations of Galois extensions condition (4) holds 
if and only if $\pi$ is an isomorphism \cite{Coh-FM90, Th. 1.2}. But $\pi$ is 
surjective by Jacobson's density theorem since $A$ is a simple left 
$A\#H$-module. Hence $\pi$ is an isomorphism if and only if $\Ker\pi=0$. If 
$\pi$ is an isomorphism then $A\#H\cong E$ is simple. On the other hand, if 
$A\#H$ is simple then $\Ker\pi=0$. Thus $(4)\Lrar(2)\Lrar(3)$.
\endproof

Condition $(1')$ has not been included as one of equivalent conditions in 
\cite{Coh-FM90, Th. 3.3}, but it appears in \cite{Coh-FM90, Cor. 3.10} stated 
for the case when the $H$-module algebra is a division ring. The proof of the 
equivalence $(1)\Lrar(2)$ given in \cite{Coh-FM90} is not quite clear since it 
is based on the equalities
$$
[A\#H/K:A^H]=[M_m(A^H):A^H]=m^2
$$
where $K$ is the kernel of the representation $\pi:A\#H\to\Endk A$, 
$\,m=[A:A^H]_r$, and $M_m(A^H)$ is the ring of $m\times m$ matrices with 
entries in $A^H$. There is a ring isomorphism $A\#H/K\cong M_m(A^H)$, but the 
natural embedding of $A^H$ in $A\#H/K$ may not correspond to the embedding of 
$A^H$ in $M_m(A^H)$ as the subring of diagonal matrices. So one has to deal 
with two different images $D$, $D'$ of the skew field $A^H$ in $M_m(A^H)$. 
Although $D\cong D'$, this does not necessarily imply that 
$$
[M_m(A^H):D]=[M_m(A^H):D'].
$$
We have given a proof which avoids the difficulty arising here.

\smallskip
If $A$ is semiprimary, then $A$ has no nontrivial $H$-stable left ideals if 
and only if $A$ has no nontrivial $H$-stable right ideals. Indeed, each of the 
two conditions imply that $A$ is $H$-simple, but then $A$ is a quasi-Frobenius 
ring by Theorem 2.8. By a standard property of quasi-Frobenius rings there is 
a bijective correspondence between left and right ideals of $A$ which assigns 
to each left ideal its right annihilator in $A$ and to each right ideal its 
left annihilator in $A$. It is easy to see that $H$-stable left ideals 
correspond to $H$-stable right ideals.

The question concerning equality of the left and right dimensions of $A$ over 
$A^H$ has been settled only recently. This question was raised in the already 
mentioned paper of Bergen, Cohen and Fischman for an $H$-module algebra which 
is a division ring, i.e. a skew field \ \cite{Ber-CF90, Question 2.4}. It was 
motivated by the classical result, due to Jacobson, that such an equality is 
indeed true in the case of group actions on skew fields.

\proclaim
Theorem 3.8 \cite{Sk-JA17}.
Suppose that $A$ is a semiprimary $H$-module algebra without nontrivial 
$H$-stable one-sided ideals. Then $\,[A\:A^H]_l=[A\:A^H]_r$.
\endproclaim

Let us outline the proof of this theorem. The desired equality holds precisely 
when the $A$-bimodule $M=\AHA$ has equal left and right ranks over $A$. By 
Lemma 3.4 $\ga$ embeds $M$ into $N=A\ot H^*$. The latter bimodule is also free 
on each side with the left and right ranks equal to $\dim H$. The conclusion 
will follow once it is shown that a direct sum of several copies of $M$ is 
isomorphic to a direct sum of several copies of $N$.

Each $A$-bimodule may be regarded as a right module over the ring
$\Ascr=A\op\ot\nobreak A$. Put $\Hscr=H\cop\ot H$ and consider $M$ and $N$ as 
left $\Hscr$-modules using the two pairs of $H$-module structures described in 
Lemma 3.2. Note that $\Hscr$ is a finite dimensional Hopf algebra and $\Ascr$ 
is an $\Hscr$-module algebra. One checks the compatibility condition which 
makes $M$ and $N$ objects of the category $\sHMA$.

Put $B=\End_\Ascr M$, i.e. $B$ is the endomorphism ring of $M$ as an 
$A$-bimodule. Then $B$ is an $\Hscr$-module algebra, and $B$ is semiprimary 
since $M$ is an $A$-bimodule of finite length. By the assumption on the 
$H$-stable one-sided ideals of $A$ there cannot exist nontrivial subbimodules 
of $M$ stable under the two $H$-module structures of Lemma 3.2. In other 
words, $M$ is a simple object of $\sHMA$. This implies $H$-semiprimeness of 
$B$, and an application of Theorem 2.4 leads further to the conclusion that 
$B$ is $H$-simple.

As a left $A$-module, and therefore as a bimodule, $N$ is generated by $1\ot 
H^*$. For each $\xi\in H^*$ there is a homomorphism of $A$-bimodules 
$\ph_\xi:M\to N$ sending $a\ot b\in\AHA\,$ to $(a\ot\xi)\,\rho(b)$ (note that 
$(1\ot\xi)\,\rho(c)=c\ot\xi$ for all $c\in A^H$), and $A\ot\xi$ is contained 
in the image of $\ph_\xi$. It follows that $N$, as an $A$-bimodule, is a 
homomorphic image of $M^d$ where $d=\dim H$. This means that the canonical map 
$$
\al:\Hom_\Ascr(M,N)\ot_BM\to N
$$
is surjective. The trickiest part is to show that $\Ker\al=0$. We refer the 
reader to \cite{Sk-JA17, Th. 3.1} for details.

Thus $N\cong F\ot_BM$ where $F=\Hom_\Ascr(M,N)$. Note that $F$ is an 
$\Hscr$-equivariant right $B$-module. By Theorem 2.1 $F^n$ is a free $B$-module 
for some $n>0$. Since $M$ and $N$ are $\Ascr$-finite, we deduce that $F$ is 
$B$-finite. Hence $F^n\cong B^r$ in $\MB$ for some integer $r>0$. This entails 
$N^n\cong F^n\ot_BM\cong B^r\ot_BM\cong M^r$ in $\sMA$, completing the proof.

In some cases the conclusion of Theorem 3.8 is almost obvious. Assume that $A$ 
is a skew field. If all $A$-bimodule composition factors of $N$ have equal 
left and right dimensions over $A$, the equality of the left and right 
dimensions over $A$ will be fulfilled for each subbimodule of $N$. This 
happens when $H$ is pointed and, more generally, when all simple $H$-comodules 
have dimension at most 2 over the ground field $\bbk$. To see this note that 
$A\ot I$ is a subbimodule of $N$ for each right ideal $I$ of $H^*$. Taking a 
composition series of $H^*$ as a right module over itself, we get a series 
of subbimodules of $A\ot H^*$ with factors isomorphic to $A\ot V$ for various 
simple right $H^*$-modules $V$ where the right action of $A$ on $A\ot V$ comes 
from $\rho$. We have $\dim V\le2$ for each $V$ by the previous assumption 
about $H$. The left and right dimensions of $A\ot V$ over $A$ are both equal 
to $\dim V$. If $\dim V=1$ then $A\ot V$ cannot contain any nontrivial 
subbimodules. If $\dim V=2$ then the $A$-bimodule $A\ot V$ is either simple or 
has exactly two composition factors, each of left and right dimension over $A$ 
equal to 1.

\medbreak
A different approach to finiteness results exploits the trace $\hatt:A\to A^H$ 
given by the action of a left integral $t\in H$ on $A$. Note that $\hatt$ is 
left and right $A^H$-linear. In particular, $\hatt$ is surjective if and only 
if $\,ta=1\,$ for some $a\in A$.

In terms of comodule structure on $A$ surjectivity of $\hatt$ is equivalent to 
the existence of a \emph{total integral} $H^*\to A$, which is a homomorphism 
of right $H^*$-comodules sending $1\in H^*$ to $1\in A$ (see \cite{Coh-F92}). 
Total integrals were introduced by Doi \cite{Doi85} for comodule algebras over 
arbitrary Hopf algebras.

\proclaim
Theorem 3.9 \cite{Mo93}.
Suppose that $A$ is right noetherian with a surjective trace $\hatt$. Then $A$ 
is a noetherian right $A^H$-module. In other words, $A^H$ is right noetherian 
and $A$ is right module-finite over $A^H$.
\endproclaim

This result is stated in \cite{Mo93, Th. 4.4.2}. To prove the theorem 
Montgomery shows that the lattice of $A^H$-submodules of any right $A$-module 
$V$ embeds into the lattice of submodules of the induced right $A\#H$-module. 
In particular, the lattice of right $A^H$-submodules of $A$ embeds into the 
lattice of right ideals of $A\#H$. Explicitly, this embedding is obtained by 
assigning to an $A^H$-submodule $U$ of $A$ the right ideal of $A\#H$ generated 
by $(U\#1)e$ where $e\in A\#H$ is an idempotent such that $A^H$ is isomorphic 
to $e(A\#H)e$.

The same argument shows that $A^H$ is right artinian if so is $A$.

The lemma below describes embeddings of lattices of submodules from our point 
of view on equivariant modules.

\proclaim
Lemma 3.10.
Let $M\in\HMA$. If the trace map $\hatt:A\to A^H$ is surjective, then\/{\rm:}

\setitemsize(iii)
\item(i)
The lattice of $A^H$-submodules of $M^H$ embeds into that of $A$-submodules of $M$.

\item(ii)
$M^H$ is an artinian $A^H$-module whenever $M$ is an artinian  $A$-module.

\item(iii)
$M^H$ is a noetherian $A^H$-module whenever $M$ is a noetherian $A$-module.

\endproclaim                      

\Proof.
We have $t(vA)=vA^H$ for any $v\in M^H$ since the map $A\to M$ such that 
$a\mapsto va$ is $H$-linear. Hence $t(UA)=U$ for each $A^H$-submodule $U$ of 
$M^H$. The assignment $U\mapsto UA$ gives therefore the desired embedding of 
lattices. Assertions (ii) and (iii) are immediate from (i).
\endproof

If $A$ is right noetherian, then each finitely generated right $A$-module is 
noetherian. Hence, if $M$ is an $A$-finite object of $\HMA$, then $M^H$ is a 
noetherian right $A^H$-module by Lemma 3.10. To deduce Theorem 3.9 from Lemma 
3.10 one needs to find an $A$-finite object $M\in\HMA$ such that $M^H\cong A$ 
as right $A^H$-modules. Recall that objects of $\HMA$ are identified with left 
$A\op\#H\cop$-modules. Let $M$ be a cyclic free $A\op\#H\cop$-module with a 
free generator $v$. Since the linear map $H\ot A\to M$ such that 
$h\ot a\mapsto h(va)$ is bijective, it follows that $M^H=tM=t(vA)$. Hence the 
assignment $a\mapsto t(va)$ defines a desired right $A^H$-linear isomorphism 
$A\to M^H$.

\setitemsize(iii)
\proclaim
Proposition 3.11.
Suppose that $H$ is semisimple and that $A$ is semiprimary and $H$-semiprime. 
Then\/{\rm:}

\item(i)
Each $H$-stable one-sided ideal of $A$ is generated by an $H$-invariant idempotent.

\item(ii)
The subring of invariants $A^H$ is semisimple artinian.

\item(iii)
$A$ is left and right module-finite over $A^H$.

\endproclaim                      

\Proof.
By Theorem 2.8 $A$ is artinian, and by Theorem 2.10 all objects of $\HMA$ are 
semisimple. In particular, the latter conclusion applies to $A$. This means 
that, whenever $I$ is an $H$-stable right ideal of $A$, there exists an 
$H$-stable right ideal $J$ such that $A=I\oplus J$. Then $I=eA$ for some 
idempotent $e\in A$. We have $e=p(1)$ where $p$ is the projection of $A$ onto 
$I$ with kernel $J$. Since $1\in A^H$ and $p$ is $H$-linear, it follows that 
$e\in A^H$. This proves assertion (i) for right ideals. Considering the 
$H\cop$-module algebra $A\op$, we get (i) for left ideals as well.

Let $U$ be any right ideal of $A^H$. Then $UA$ is an $H$-stable right ideal of 
$A$. By (i) $UA$ is generated by an idempotent $e\in A^H$. By Lemma 3.10 applied 
to $M=A$ the lattice of right ideals of $A^H$ embeds into that of right ideals 
of $A$. Since the right ideal $eA^H$ of $A^H$ has the same extension to $A$ as 
$U$, we deduce that $U=eA^H$. Thus each right ideal of $A^H$ is generated by 
an idempotent. This yields (ii). Finally, (iii) follows from Theorem 3.9.
\endproof

We have given a self-contained proof. It has been known for a long time that 
$A^H$ is semisimple artinian whenever so is $A\#H$ \cite{Coh-FM90, Th. 3.13}. 
That $A\#H$ is semisimple artinian, under the hypothesis of Proposition 3.11, 
has been established in \cite{Sk-Oy06} (see Theorem 2.10).

\section
4. Localization at invariants and a Bergman-Isaacs type theorem

In \cite{Sk15} it was shown that for a semisimple Hopf algebra $H$ all right 
noetherian \hbox{$H$-module} algebras have, loosely speaking, ``sufficiently 
many" $H$-invariant elements. The proofs of these results referred to a few 
statements from \cite{Sk-Oy06} which served as intermediate steps in the 
process of verifying the existence of artinian classical quotient rings. But 
in fact only the final conclusion from \cite{Sk-Oy06} is needed, and therefore 
those results of \cite{Sk15} are valid for a larger class of $H$-module 
algebras. This will be explained below.

Suppose that an $H$-semiprime $H$-module algebra $A$ has a right artinian 
classical right quotient ring $Q$. By Theorem 2.14 the $H$-module structure 
extends to $Q$. It is clear then that $Q$ has to be $H$-semiprime since 
$J\cap A\ne0$ for each nonzero right ideal $J$ of $Q$. Therefore all results 
concerning artinian $H$-semiprime algebras apply to $Q$.

A right ideal $I$ of a ring $R$ is called \emph{essential} if $I$ has nonzero 
intersection with each nonzero right ideal of $R$.

\setitemsize(a)
\proclaim
Lemma 4.1.
Suppose that $A$ is $H$-semiprime and $A$ has a right artinian classical right 
quotient ring $Q$. For a right ideal $I$ of $A$ denote by $I_H$ the largest 
$H$-stable right ideal of $A$ contained in $I$. The following conditions 
are equivalent{\rm:}

\item(a)
$I_H$ is an essential right ideal of $A$,

\item(b)
$I$ contains a regular element of $A$.

\endproclaim

\Proof.
Suppose that $I_H$ is an essential right ideal of $A$. Then $I_HQ$ is an 
essential right ideal of $Q$. But $I_HQ$ is $H$-stable, and therefore this right 
ideal is generated by an idempotent $e\in Q$ according to Corollary 2.7. Since 
$eQ\cap(1-e)Q=0$, we must have $e=1$, i.e $I_HQ=Q$. Then $1=as^{-1}$ for some 
$a\in I_H$ and a regular element $s\in A$. So $s=a\in I_H$, which proves that 
(a)${}\Rar{}$(b).

Suppose now that $I$ contains a regular element $s$ of $A$. Note that
$$
I_H=\{a\in A\mid Ha\sbs I\}=\{a\in A\mid\rho(a)\in I\ot H^*\}.
$$
We have to prove that $I_H$ is an essential right ideal of $A$. Suppose that 
$I_H\cap bA=0$ for some $b\in A$, $\,b\ne0$. Then $(I\ot H^*)\cap\rho(bA)=0$. 
In particular,
$$
(sA\ot H^*)\cap\rho(bA)=0.
$$
Since $s\ot1$ is a regular element of the ring $A\ot H^*$, it follows that the sum 
$$
\sum_{n=0}^\infty\,(s^n\ot1)\,\rho(bA)
$$
is direct, and each summand is a nonzero right $\rho(A)$-submodule of $A\ot H^*$. 
Consider now $F=A\ot H^*$ as a right $A$-module with respect to the action of $A$ 
given by right multiplications by the elements $\rho(a)$, $\,a\in A$. We know 
that this $A$-module is free of rank equal to the dimension of $H$. Hence 
$F\ot_AQ$ is a finitely generated right $Q$-module containing an infinite 
direct sum of nonzero submodules. However, this is impossible since $Q$ is 
right artinian. Thus $I_H\cap bA\ne0$ whenever $b\ne0$, and so (b)$\Rar$(a).  
\endproof

\setitemsize(ii)
\proclaim
Lemma 4.2.
Suppose that $H$ is semisimple, $A$ is $H$-semiprime, and $A$ has a right 
artinian classical right quotient ring $Q$. Let $I$ be an $H$-stable right 
ideal of $A$, and let $\,I^H=I\cap A^H$.

\item(i)
If $I^H=0$ then $I=0$.

\item(ii)
If $I$ is an essential right ideal of $A$ then $Q^HI^H=Q^H$ and $I^HQ^H=Q^H$.

\endproclaim

\Proof.
By Proposition 3.11 $Q^H$ is semisimple artinian and $Q$ has finite length as 
either left or right $Q^H$-module.

Consider the $(Q^H\!\!,A)$-subbimodule $Q^HI$ of $Q$. For each regular element 
$s$ of $A$ the $Q^H$-submodule $Q^HIs$ is isomorphic to $Q^HI$. Hence these 
two $Q^H$-modules have the same length. Since $Q^HIs\sbs Q^HI$, we get 
$Q^HIs=Q^HI$, and therefore $Q^HIs^{-1}=Q^HI$.

It follows that $Q^HI$ is a right ideal of $Q$. Since it is $H$-stable, 
Proposition 3.11 yields $Q^HI=eQ$ for some $e\in Q^H$. Let $t\in H$ be an 
integral. The action of $t$ on $Q$ commutes with the left and right 
multiplications by $H$-invariant elements. Hence
$$
e\in eQ^H=t(eQ)=t(Q^HI)=Q^HI^H.
$$
If $I^H=0$, the above inclusion entails $e=0$, i.e. $I=0$. This proves (i).

Suppose that $I$ is an essential right ideal of $A$. By Lemma 4.1 $I=I_H$ 
contains a regular element of $A$, so that $IQ=\nobreak Q$. Since $Q^HI$ is a 
right ideal of $Q$ containing $I$, we get $Q^HI=Q$ as well. The previous 
argument with $e=1$ shows that $1\in Q^HI^H$. Thus $Q^H=Q^HI^H$.

Consider now the $H$-stable right ideal $I^HQ$ of $Q$. By Theorem 2.10 there 
exists an $H$-stable right ideal $J$ of $Q$ such that $Q=I^HQ\oplus J$. Then 
$I\cap J$ is an $H$-stable right ideal of $A$ with
$$
(I\cap J)^H=I^H\cap J=0.
$$
As we have proved already in part (i) this entails $I\cap J=0$. Since $IQ=Q$, 
any element $y\in J$ can be written as $y=as^{-1}$ where $a\in I$ and $s$ is a 
regular element of $A$; then $a\in I\cap J$, so that $a=0$ and $y=0$. 
Therefore $J=0$ and $Q=I^HQ$. Hence
$$
Q^H=tQ=t(I^HQ)=I^HQ^H,
$$
and we are done.
\endproof

\setitemsize(iii)
\proclaim
Theorem 4.3.
Suppose that $H$ is semisimple, $A$ is $H$-semiprime, and $A$ has a right 
artinian classical right quotient ring $Q$. Denote by $\Si$ the set of regular 
elements of $A^H$ and by $\E$ the set of right ideals of $A$ which satisfy 
the equivalent conditions {\rm(a)} and {\rm(b)} of Lemma\/ {\rm4.1}. 
Then\/{\rm:}

\item(i)
The algebra $A^H$ is semiprime right Goldie.

\item(ii)
$\Si$ is a right Ore subset of regular elements of $A$.

\item(iii)
$Q$ is canonically isomorphic with the right localization of $A$ at $\Si$.

\item(iv)
The classical right quotient ring of $A^H$ is isomorphic with $Q^H$.

\item(v)
$I\cap\Si\ne\noth$ for each right ideal $I\in\E$.

\endproclaim

\Proof.
Put $I^H=I\cap A^H$ and $\F=\{I^H\mid I\in\E\}$. If $I,J\in\E$ then 
$I\cap J\in\E$, whence $I^H\cap J^H=(I\cap J)^H\in\F$. Therefore $\F$ is a 
filter of right ideals of $A^H$. If $I\in\E$ and $a\in A^H$, then the right 
ideal $I_a=\{x\in A\mid ax\in I_H\}$ is essential and $H$-stable; so 
$I_a\in\E$ and $I_a^H\in\F$. Also, $aI_a^H\sbs I^H$ since $aI_a\sbs I_H\sbs I$. 
This shows that $\F$ is a topologizing filter.

Moreover, $\F$ satisfies conditions (a) and (b) of Proposition 2.13 
with $A$ replaced by $A^H$ and $Q$ replaced by $Q^H$. Indeed, 
if $I\in\E$ then $I^H$ has zero left and right annihilators in $Q$ since 
$I^HQ$ and $QI^H$ contain the unity 1 by Lemma 4.2. If $q\in Q^H$ there exists 
a regular element $s$ of $A$ such that $qs\in A$. Put $K=sA$. Then $K\in\E$ by 
Lemma 4.1 and $qK\sbs A$. Hence $K^H\in\F$ and $qK^H\sbs A^H$.

It has been observed already that the ring $Q^H$ is semisimple artinian. Now 
(iv) follows from Proposition 2.13, and (i) is its consequence since a ring $R$ is 
semiprime right Goldie if and only if $R$ has a semisimple artinian classical 
right quotient ring.

Given $I\in\E$, the equality $I^HQ^H=Q^H$ of Lemma 4.2 means that 
$I^H\cap\Si\ne\noth$, which amounts to (v). For any $q\in Q$ the set 
$$
I=\{x\in A\mid qx\in A\}
$$
is a right ideal of $A$ containing a regular element of $A$. By Lemma 4.1 
$I\in\E$. Since $qI\sbs A$, assertion (v) shows that $qs\in A$ for some 
$s\in\Si$.

All elements of $\Si$ are invertible in $Q^H$ and therefore in $Q$. Hence all 
elements of $\Si$ are regular in $A$. Since each element of $Q$ can be written 
in the form $as^{-1}$ for some $a\in A$ and $s\in\Si$, (ii) and (iii) are 
immediate (see \cite{McC-R, 2.2.4}).
\endproof

In all corollaries below we continue to assume that $H$ is semisimple.

\proclaim
Corollary 4.4.
Let $A$ be as in Theorem {\rm4.3}. If $s$ is any regular element of $A$,
then the right ideal $sA$ contains an $H$-invariant regular element of $A$.
\endproclaim

\Proof.
Since $sA\in\E$ by Lemma 4.1, we have $sA\cap\Si\ne\noth$ by Theorem 4.3.
\endproof

\proclaim
Corollary 4.5.
Let $A$ and $Q$ be as in Theorem {\rm4.3}, and let $I$ be any $H$-stable right 
ideal of $A$. There exists $x\in I^H$ such that $IQ=xQ$.
\endproclaim

\Proof.
Since $IQ$ is an $H$-stable right ideal of $Q$, we have $IQ=eQ$ for some 
$e\in Q^H$ by Proposition 3.11. Now $J=\{a\in A\mid ea\in I\}$ is an 
$H$-stable right ideal of $A$ containing a regular element of $A$. Hence 
$J\cap\Si\ne\noth$ by Corollary 4.4. Pick any $s\in J\cap\Si$ and put $x=es$. 
Then $x\in I\cap Q^H=I^H$. Since $s$ is invertible in $Q$, we get $eQ=xQ$.  
\endproof

\setitemsize(a)
\proclaim
Corollary 4.6.
All conclusions of Theorem {\rm4.3,} as well as Corollaries {\rm4.4, 4.5,} 
hold in each of the following three cases\/{\rm:}

\item(a)
$A$ is semiprime right Goldie,

\item(b)
$A$ is right noetherian and $H$-semiprime,

\item(c)
$A$ is PI and $H$-semiprime with finitely many minimal $H$-prime ideals.

\endproclaim

\Proof.
A right artinian classical right quotient ring $Q$ exists in case (a) by the 
Goldie theorem, in the other cases by Theorems 2.11, 2.12.
\endproof

Under the assumption that $A\#H$ is semiprime, a short argument given by Bergen 
and Montgomery \cite{Ber-M86, Prop. 2.4} shows that $A^H$ is semiprime, and 
that $\hatt(I)\ne0$ where $\hatt:A\to A^H$ is the trace map (in particular, 
$I^H\ne0$) for each nonzero \hbox{$H$-stable} one-sided ideal $I$ of $A$. From 
this it was further deduced in \cite{Ber-M86, Lemma 3.4} that, among other 
things, regular elements of $A^H$ are regular in $A$, and that $A^H$ is Goldie 
when so is $A$. If $A,Q,H$ are as in Theorem 4.3, then $Q\#H$ is semisimple 
artinian by Theorem 2.10; since $Q\#H$ is a classical right quotient ring of 
$A\#H$, it follows that $A\#H$ is semiprime. This fact was not known at the 
time when \cite{Ber-M86} was written.

Several deeper results from \cite{Ber-M86} use the assumption that $A\#H$ is 
not only semiprime, but has the \emph{ideal intersection property} (\IIP for 
short) which means that each nonzero ideal of $A\#H$ has nonzero intersection 
with $A$. In fact, in the presence of \IIP the ring $A\#H$ is semiprime if and 
only if $A$ is $H$-semiprime. The \IIP is satisfied for $X$-outer group actions 
on semiprime rings and for $X$-outer actions of Lie algebras on prime rings. 
However, it seems that there are no approaches to analogs of such results for 
actions of arbitrary finite dimensional or even semisimple Hopf algebras.

It was asked in \cite{Ber-M86} whether $Q(A)^H=Q(A^H)$ when $A\#H$ is semiprime 
with \IIP. Part (iv) of Theorem 4.3 answers this question, imposing reasonable 
conditions on $A$ and $H$, but not assuming the \IIP. In the case of a finite 
group $G$ acting on a semiprime ring $R$ without additive $|G|$-torsion the 
fact that $R^G$ is right Goldie if and only if $R$ is right Goldie and the 
equality $Q(R)^G=Q(R^G)$ were proved by Kharchenko \cite{Kh74}; it was also 
observed by Montgomery \cite{Mo78} that $Q(R)$ is the localization of $R$ at 
the Ore set of regular $G$-invariant elements. Analogs of these results for 
group graded rings are due to Cohen and Rowen \cite{Coh-R83}.

\medbreak
There is a slightly weaker version of Lemma 4.2 for $H$-stable left ideals. 
The equality $Q^HI^H=Q^H$ in (ii) cannot be proved unless $Q$ is a two-sided 
quotient ring.

\setitemsize(ii)
\proclaim
Lemma 4.7.
Suppose that $H$ is semisimple, $A$ is $H$-semiprime, and $A$ has a right 
artinian classical right quotient ring $Q$. Let $I$ be an $H$-stable left 
ideal of $A$, and let $\,I^H=I\cap A^H$.

\item(i)
If $I^H=0$ then $I=0$.

\item(ii)
If $QI=Q$ then $I^HQ^H=Q^H$.

\endproclaim

\Proof.
We repeat the steps in the proof of Lemma 4.2 using the two-sided properties 
of $Q$. First, the $(A,Q^H)$-subbimodule $IQ^H$ is a left ideal of $Q$ since $Q$ 
has finite length as a right $Q^H$-module. Next, $IQ^H=Qe$ for some $e\in Q^H$ 
by Proposition 3.11. Applying the integral $t$, we deduce that $e\in I^HQ^H$. 

If $I^H=0$, then $e=0$, whence $I=0$. If $QI=Q$, then $IQ^H=Q$, and therefore 
$1\in I^HQ^H$.
\endproof

\proclaim
Theorem 4.8.
Denote by $N$ the prime radical of $A$ and by $N_H$ the largest $H$-stable 
ideal of $A$ contained in $N$. Suppose that $H$ is semisimple, $N$ is 
nilpotent, and $A/N_H$ has a right artinian classical right quotient ring. If 
$I$ is any $H$-stable one-sided ideal of $A$ such that $I^H$ is nilpotent, 
then $I$ is nilpotent.
\endproclaim

\Proof.
Denote by $\pi$ the canonical surjective homomorphism of $H$-module algebras  
$A\to A/N_H$. Then $\pi(I)$ is an $H$-stable one-sided ideal of $A/N_H$ and 
$\pi$ maps $I^H$ onto $\pi(I)^H\!$ since $H$ is semisimple. Hence $\pi(I)^H$ 
is nilpotent. The factor algebra $A/N_H$ is $H$-semiprime since $N_H$ is an 
$H$-semiprime ideal of $A$.  Hence its subring of invariants $(A/N_H)^H$ is 
semiprime by Theorem 4.3, which entails $\pi(I)^H=0$. An application of Lemmas 
4.2, 4.7 yields $\pi(I)=0$, i.e. $I\sbs N$. So $I$ is nilpotent.
\endproof

\setitemsize(a)
\proclaim
Corollary 4.9.
Suppose that $H$ is semisimple. The conclusion of Theorem {\rm4.8} holds in 
each of the following three cases\/{\rm:}

\item(a)
$A$ is left noetherian,

\item(b)
$A$ is right noetherian,

\item(c)
$A$ is finitely generated and PI.

\endproclaim

\Proof.
In all cases $N$ is known to be nilpotent. In cases (b) and (c) the $H$-semiprime 
factor algebra $A/N_H$ has a right artinian classical right quotient ring by 
Theorems 2.11, 2.12. In case (a) we apply Theorem 4.8 to the right noetherian 
$H\cop$-module algebra $A\op$.
\endproof

Let $R$ be a nonunital ring and $G$ a finite group of its automorphisms 
such that $R$ has no additive $|G|$-torsion. Consider the trace map 
$$
\hatt:R\to R^G,\qquad\hatt(a)=\tsum\nolimits_{g\in G}\,ga.
$$
A classical result of Bergman and Isaacs \cite{Ber-I73} says that $R$ is 
nilpotent if $\hatt(R)$ is nilpotent. Moreover, if $\hatt(R)=0$, then the 
nilpotency index of $R$ is bounded by a number which depends only on the order 
$|G|$ of the group $G$, but not on the ring $R$. An easy consequence of this 
result is that $R^G$ is semiprime when so is $R$.

A similar result for group graded rings proved by Cohen and Rowen \cite{Coh-R83} 
is even simpler: if $R=\bigoplus_{g\in G}R_g$ is a nonunital ring graded by a 
finite group $G$, and if $R_1^d=0$ for some integer $d>0$, then $R^{\,d|G|}=0$. 
In fact one needs only a grading with finite support, while the group $G$ may 
be infinite.

Theorem 4.8 is different not only in that the conclusion is stated for one-sided 
ideals of $H$-module algebras satisfying certain conditions, but also because 
it relies heavily on the $H$-semiprime case. The nilpotency index of $I$ can 
be bounded only by the nilpotency index of the ideal $N_H$, even when $I^H=0$.

Bahturin and Linchenko \cite{Bah-L98} investigated conditions under which one 
can conclude that $A$ is PI, knowing that $A^H$ is PI. They showed that, for a 
fixed finite dimensional Hopf algebra $H$, in order that each $H$-module algebra 
$A$ be PI whenever $A^H$ is PI it is necessary and sufficient that there exist 
a natural number $n$ such that $A^n=0$ for each nonunital $H$-module algebra 
$A$ with $A^H\cdot A^H=0$, and this can happen only if $H$ is semisimple. 
Several other equivalent conditions are given in \cite{Bah-L98}. This work of 
Bahturin and Linchenko elucidates the need for a more precise analog of the 
Bergman-Isaacs result for Hopf algebra actions.

\section
5. Hopf actions on commutative algebras

Throughout this section we assume that $A$ is a \emph{commutative} $H$-module 
algebra. First we are going to recall the algebraic interpretation of the 
classical result on quotients of affine schemes by actions of finite group 
schemes.

Given an associative algebra $U$ over a commutative ring $R$ such that $U$ is 
free of finite rank as an $R$-module the \emph{norm}\/ $\Nm_{U/R}(u)\in R$ of an 
element $u\in U$ is defined as the determinant of the operator $L_u\in\End_RU$ 
of the left multiplication by $u$ in $U$. Considering the polynomial ring 
$U[t]$, where $t$ is an indeterminate, as an algebra over $R[t]$ we get also 
the \emph{characteristic polynomial}
$$
P_{U/R}(u,t)=\Nm_{U[t]/R[t]}(t-u)=\det(t\cdot\Id-L_u)\in R[t].
$$
In particular, $(-1)^r\Nm_{U/R}(u)$ where $r=\rank_RU$ is the coefficient of 
$t^0$ in this polynomial. Passing to localizations of the base ring $R$ these 
definitions extend to the case where $U$ is not free as an $R$-module, but 
only projective of finite constant rank.

By the Cayley-Hamilton theorem $P_{U/R}(u,L_u)=0$ in $\End_RU$. Applying this 
operator to the identity element $1\in U$, we get $P_{U/R}(u,u)=0$ in $U$, 
which is a relation of integral dependence of the element $u$ over the ring 
$R$. The integral dependence of $U$ over $R$ is merely a consequence of 
module-finiteness. What is important for application to invariants is the fact 
that the characteristic polynomials enjoy several nice properties. In 
particular, they are functorial in the sense that, given a homomorphism of 
commutative rings $\ze:R\to R'$, we have
$$
P_{R'\ot_RU/R'}(1\ot u,t)=\ze^tP_{U/R}(u,t)
$$
where $\ze^t:R[t]\to R'[t]$ is the homomorphism extending $\ze$ and sending 
$t$ to $t$.

Since $A$ is commutative, the map
$$
\io:A\to A\ot H^*,\qquad a\mapsto a\ot1,
$$
is an isomorphism of $A$ onto a central subalgebra of $A\ot H^*$. So we may 
regard $A\ot H^*$ as an $A$-algebra via $\io$. Clearly this algebra is free of 
rank $d=\dim H$ as an $A$-module. In this way the polynomial $P_{A\ot H^*/A}(u,t)$ 
is defined for each $u\in A\ot H^*$. Making use of the comodule structure 
$\rho:A\to A\ot H^*$, we get the polynomial
$$
P_{A\ot H^*/A}\bigl(\rho(a),t\bigr)\in A[t]
$$
for $a\in A$. Suppose that 
$$
P_{A\ot H^*/A}\bigl(\rho(a),t\bigr)=t^d+\tsum_{i=0}^{d-1}c_it^i\quad
\hbox{where $c_0,\ldots,c_{d-1}\in A$}.
$$
Then $\rho(a)^d+\sum_{i=0}^{d-1}(c_i\ot1)\rho(a)^i=0$ in $A\ot H^*$. Applying 
to the left hand side of this equality the algebra homomorphism 
$\id\ot\ep:A\ot H^*\to A$, we get 
$$
a^d+\tsum_{i=0}^{d-1}c_ia^i=0.
$$
If $c_i\in A^H$ for all $i$, the relation above shows that $a$ is integral 
over $A^H$. On this observation the classical argument reproduced in the 
following theorem is based:

\proclaim
Theorem 5.1.
Suppose that $H$ is cocommutative. Then for each $a\in A$ the characteristic 
polynomial $\,P_{A\ot H^*/A}\bigl(\rho(a),t\bigr)\,$ has all coefficients in 
$A^H$. In particular, $A$ is integral over $A^H$.
\endproclaim

\Proof.
By the condition imposed on $H$ the dual Hopf algebra $H^*$ is commutative, 
whence $A\ot H^*$ is commutative as well. Since $A^H$ is the equalizer of 
the two algebra homomorphisms $\,\io,\,\rho:A\to A\ot H^*$, we have to show that
$$
\io^tP_{A\ot H^*/A}\bigl(\rho(a),t\bigr)=
\rho^tP_{A\ot H^*/A}\bigl(\rho(a),t\bigr).\eqno(*)
$$
Note that the commutative diagrams
$$
\diagram{
A&\hidewidth\lmapr6{\io}\hidewidth&A\ot H^*&
&A&\hidewidth\lmapr6{\io}\hidewidth&A\ot H^*&\cr
\noalign{\smallskip}
\lmapd{16pt}{\io}{}&&\lmapd{16pt}{}{\io\ot\id}&&
\lmapd{16pt}{\rho}{}&&\lmapd{16pt}{}{\rho\ot\id}\cr
\noalign{\smallskip}
A\ot H^*&\lmapr2{\io'}&A\ot H^*\ot H^*&\hskip4em
&A\ot H^*&\lmapr2{\io'}&A\ot H^*\ot H^*&\cr
}
$$
where $\io'(x)=x\ot1$ for $x\in A\ot H^*$ are cocartesian in the category of 
commutative $A$-algebras in the sense that each diagram makes $A\ot H^*\ot H^*$ 
the tensor product of two $A$-algebras given by the respective homomorphisms 
$A\to A\ot H^*$. Hence $(*)$ can be rewritten as
$$
P_{A\ot H^*\ot H^*/A\ot H^*}\bigl((\io\ot\id)\rho(a),t\bigr)=
P_{A\ot H^*\ot H^*/A\ot H^*}\bigl((\rho\ot\id)\rho(a),t\bigr)\eqno(**)
$$
by functoriality of the characteristic polynomials. Here $A\ot H^*\ot H^*$ 
is viewed as an $A\ot H^*$-algebra by means of $\io'$. Next, there is an 
automorphism $\ph$ of the algebra $A\ot H^*\ot H^*$ defined by the rule 
$$
\ph(x\ot\xi)=(x\ot1)(1\ot\De\xi)\quad\hbox{for $x\in A\ot H^*$ and 
$\xi\in H^*$}.
$$
Since $\ph$ acts as the identity on $A\ot H^*\ot1$, we have
$$
P_{A\ot H^*\ot H^*/A\ot H^*}(y,t)=
P_{A\ot H^*\ot H^*/A\ot H^*}\bigl(\ph(y),t\bigr)
$$
for all $y\in A\ot H^*\ot H^*$. If $y=(\io\ot\id)\rho(a)$, then
$
\ph(y)=(\id\ot\De)\rho(a)=(\rho\ot\id)\rho(a).
$
Hence $(**)$ follows.
\endproof

By passage from $A$ to $A[t]$ the conclusion of Theorem 5.1 can be deduced 
from the fact that
$$
\Nm_{A\ot H^*/A}\bigl(\rho(a),t\bigr)\in A^H\quad\hbox{for all $a\in A$}.
$$
The arguments showing this inclusion in the proof of Theorem 5.1 above are 
tautological, up to different notation and the use of right comodule 
structures instead of left ones, to those in Mumford's book on abelian 
varieties \cite{Mum, Ch. III, section 12}. Demazure and Gabriel describe 
quotients by actions of finite group schemes in \cite{De-G, Ch. III, \S2, 
Cor. 6.1} as a special case of a more general result on quotients of groupoid 
schemes \cite{De-G, Ch. III, \S2, Th. 3.2}.

Thus Theorem 5.1 is a very old result. Somehow it had not been well-known to 
Hopf algebra theorists for some time in the past. Integrality over invariants 
for commutative comodule algebras over commutative Hopf algebras was 
rediscovered by Ferrer Santos in \cite{Fer94}. In the language of module 
algebras that approach was reformulated by Montgomery \cite{Mo93, \S4.2}. It 
makes use of the characteristic polynomials of endomorphisms of equivariant 
$A$-modules.

In \cite{Mo93} Montgomery raised the question as to whether $A$ is always 
integral over $A^H$ in the case of an arbitrary finite dimensional Hopf 
algebra $H$. For pointed Hopf algebras this question was answered shortly 
afterwards in the affirmative by Artamonov \cite{Art96} when $A$ is a domain 
and without any restrictions on $A$ by Totok \cite{Tot97} and Zhu \cite{Zhu96} 
when $\chr\bbk>0$. Both \cite{Tot97} and \cite{Zhu96} provided counterexamples 
to integrality in characteristic 0. Zhu also proved that $A$ is integral over 
$A^H$ when $H$ is involutory, i.e. $S^2=\id$, and $\chr\bbk$ does not divide the 
dimension of $H$. At that time it remained open what is actually needed for 
integrality to hold.

The characteristic polynomials have reappeared in a later development:

\proclaim
Theorem 5.2 \cite{Sk04}.
If $A$ is $H$-semiprime or, more generally, if there exists a homomorphism of 
commutative $H$-module algebras $\ph:A'\to A$ such that $A=\ph(A')A^H$ and 
$A'$ is $H$-semiprime then for each $a\in A$ the polynomial 
$\,P_{A\ot H^*/A}\bigl(\rho(a),t\bigr)\,$ has all coefficients in $A^H$. In 
particular, $A$ is integral over $A^H$.
\endproclaim

The case when $A$ is $H$-semiprime, which is the main step here, has been 
subsumed in a recent work of Eryashkin \cite{Er-a16} on invariants of 
$H$-module PI algebras. These results will be discussed in section 7. The 
original proof of Theorem 5.2 had common elements with the proof of Theorem 
7.5, but it didn't use the Martindale quotient rings.

If $A$ contains nonzero $H$-stable nilpotent ideals, then integrality over 
invariants may well be lost by the already mentioned examples of Totok and Zhu. 
There are still two important cases when the $H$-semiprimeness is not needed:

\proclaim
Corollary 5.3.
$A$ is integral over $A^H$ in each of the following two cases{\rm:}

\item(a)
the trace map $A\to A^H$ is surjective,

\item(b)
$\chr\bbk=p>0$.

\endproclaim

\Proof.
Let $N$ be the largest $H$-stable ideal of $A$ contained in the nil radical of 
$A$. Since $B=A/N$ is $H$-semiprime, Theorem 5.2 shows that $B$ is integral 
over $B^H$. Let $\pi:A\to B$ be the canonical map.

In case (a) $\pi(A^H)=B^H$. Therefore for each $a\in A$ there exists a 
polynomial $f\in A^H[t]$ with the leading coefficient 1 such that $f(a)\in N$. 
Then $f(a)^n=0$ for some integer $n>0$ since $N$ is nil. Hence $a$ is integral 
over $A^H$.

Suppose that $\chr\bbk=p>0$. Put $A'=\pi^{-1}(B^H)$. As in case (a) it is 
checked that $A$ is integral over $A'$. We claim that for each $c\in A'$ there 
exists $n>0$ such that $c^{p^n}\in A^H$. Indeed, $\rho(c)-c\ot1\in N\ot H^*$ 
since $\pi(c)\in B^H$. It follows that $\rho(c)-c\ot1$ is nilpotent. Hence
$$
\rho(c^{p^n})-c^{p^n}\ot1=(\rho(c)-c\ot1)^{p^n}=0
$$
for sufficiently large $n$, but this means that $c^{p^n}\in A^H$. Thus $A'$ is 
 integral over $A^H$, and the final conclusion follows from transitivity of 
integrality.
\endproof

In \cite{Zhu96} Zhu conjectured that $A$ is integral over $A^H$ whenever $H$ 
is involutory. When $\chr\bbk=0$ it is known that $H$ is involutory if and only 
if $H$ is semisimple. In this case the trace $A\to A^H$ is surjective. Thus 
Zhu's conjecture follows from Corollary 5.3. However, when $\chr\bbk>0$ the 
question of integrality does not depend on any condition on $H$.

As observed by Kalniuk and Tyc \cite{Kal-T08} the fact that in positive 
characteristic each commutative $H$-module algebra is integral over the 
invariants implies a property of $H$ similar to the geometric reductivity 
known in the theory of algebraic groups. This property was considered in 
\cite{Kal-T08} for a not necessarily finite dimensional Hopf algebra $H$, and 
its main consequence is that, whenever $A$ is a finitely generated commutative 
$H$-module algebra on which the action of $H$ is locally finite, the algebra 
of invariants $A^H$ is finitely generated. When $\chr\bbk>0$ each finite 
dimensional Hopf algebra is geometrically reductive in this sense 
\cite{Kal-T08, Th. 4}. This result can be reformulated as follows:

\proclaim
Theorem 5.4 \cite{Kal-T08}.
Suppose that $\chr\bbk=p>0$. If $\ph:A\to B$ is a surjective homomorphism 
of commutative $H$-module algebras and $b\in B^H$, then $b^n\in\ph(A^H)$ 
for some integer $n>0$.
\endproclaim

\Proof.
Consider first the case when $A$ and $B$ are graded, $\ph$ respects the 
grading, and $b$ is homogeneous of degree 1. Let $B'$ be the subalgebra of $B$ 
generated by $b$, and let $A'=\ph^{-1}(B')$. Since $A'$ is integral over 
$A'^H$ by Corollary 5.3, $B'=\ph(A')$ is integral over $\ph(A'^H)$. But  
$\ph(A'^H)$ is a graded subalgebra of $B'$. If $\ph(A'^H)=\bbk$, then $b$ is 
integral over $\bbk$, in which case $b$ has to be nilpotent. Otherwise 
$\ph(A'^H)$ contains a homogeneous element of positive degree. But each 
homogeneous component of $B'$ is spanned by a power of $b$. Hence 
$b^n\in\ph(A^H)$ for some $n$ anyway.

In the general case let $\ph^t:A[t]\to B[t]$ be the extension of $\ph$ to 
polynomial rings and extend the action of $H$ to $A[t]$ and $B[t]$ by making 
$t$ invariant. Then $bt\in B[t]$ is a homogeneous $H$-invariant element of 
degree 1, and so we are in the situation of the previous case.
\endproof

Integrality of $A$ over $A^H$ implies several well-known nice properties of 
the ring extension $A^H\sbs A$. In particular, the canonical map $\Spec 
A\to\Spec A^H$ between the prime spectra is surjective, closed, and satisfies 
the going-up. However, for deeper conclusions integrality alone is not 
sufficient, and the characteristic polynomials come into play in an essential 
way. One application is this:

\proclaim
Theorem 5.5 \cite{Sk04}.
Suppose that for each $a\in A$ the characteristic polynomial 
$\,P_{A\ot H^*/A}\bigl(\rho(a),t\bigr)\,$ has all coefficients in $A^H$. Then 
the map\/ $\Spec A\to\Spec A^H$ is open, has finite fibers, and satisfies 
the going-down property.
\endproclaim

Theorem 5.5 and its proof generalize the classical results describing 
properties of the quotient morphism $X\to X/G$ where $X$ is an affine scheme 
and $X/G$ is its quotient by an action of a finite group scheme.

There are further applications of the technique used in the study of group 
scheme actions. For each $\p\in\Spec A$ denote by $k(\p)$ the residue field 
of the local ring $A_\p$. Let $\al_\p:A\to k(\p)$ be the canonical ring 
homomorphism. The composite
$$
\de_\p:A\lmapr2\rho A\ot H^*\lmapr6{\al_\p\ot\id}k(\p)\ot H^*
$$
is a homomorphism of $H$-module algebras, assuming that $H$ acts trivially on 
$k(\p)$ and by the left hits $\rhu$ on $H^*$. Hence
$$
O(\p)=\bigl(k(\p)\ot1\bigr)\cdot\de_\p(A)
$$
is a commutative right coideal subalgebra of the Hopf algebra $k(\p)\ot H^*$ 
over the field $k(\p)$. In \cite{Sk04} $O(\p)$ was called the \emph{orbital 
subalgebra} associated with $\p$.

When $H$ is cocommutative and $G$ is the finite group scheme representable by 
the commutative Hopf algebra $H^*$, the algebra $O(\p)$ represents the 
scheme-theoretic $G$-orbit of $\p$ which is a closed subscheme in the affine 
scheme $\Spec(k(\p)\ot A)$.

\proclaim
Theorem 5.6 \cite{Sk04}.
Suppose that $A$ is $H$-semiprime and the function\/ $\p\mapsto\dim O(\p)$ is 
locally constant on the whole $\Spec A$. Then $A$ is a finitely generated 
projective $A^H$-module whose rank at a prime\/ $\q\in\Spec A^H$ is equal to 
$\dim O(\p)$ where $\p$ is any prime ideal of $A$ lying above $\q$.

Also, the assignment $I\mapsto I\cap A^H$ establishes a bijection between the 
$H$-stable ideals of $A$ and all ideals of $A^H$. The inverse correspondence 
is $J\mapsto JA$.
\endproclaim

We will explain briefly the main ideas used in the proof. Given some elements 
$a_1,\ldots,a_n\in A$, the set $U$ of those prime ideals $\p$ of $A$ for which 
$\de_\p(a_1),\ldots,\de_\p(a_n)$ form a basis of $O(\p)$ over $k(\p)$ is open 
in $\Spec A$. One can also check that, whenever $\p$ and $\p'$ are two prime 
ideals of $A$ with $\p\cap A^H=\p'\cap A^H$, one has $\p\in U$ if and only if 
$\p'\in U$. Then, passing to the localizations $A[s^{-1}]$ of $A$ at suitable 
elements $s\in A^H$, one may assume that there exist $a_1,\ldots,a_n\in A$ 
such that $\de_\p(a_1),\ldots,\de_\p(a_n)$ are a basis of $O(\p)$ over 
$k(\p)$ for each $\p\in\Spec A$.

The technically complicated part of the proof is to show that the previous 
assumption implies that $\rho(a_1),\ldots,\rho(a_n)$ form a basis of 
$(A\ot1)\rho(A)\sbs A\ot H^*$ over $A$ with respect to the left module 
structure; once this has been done, the freeness of $A$ over $A^H$ follows from 
Lemma 3.3. Note, however, that, when $A$ is reduced (equivalently, semiprime), 
there is a general ring-theoretic fact which states that a submodule $M$ of a 
finite rank free $A$-module $F$ is freely generated by elements 
$v_1,\ldots,v_n\in M$ provided that for each $\p\in\Spec A$, the image of 
$k(\p)\ot M$ in $k(\p)\ot F$ has a basis over $k(\p)$ consisting of $1\ot 
v_1,\ldots,1\ot v_n$; with
$$
M=(A\ot1)\rho(A),\qquad F=A\ot H^*,\qquad v_i=\rho(a_i)
$$
the desired conclusion is immediate. Since $A$ is assumed to be only 
$H$-semiprime, one has to overcome several difficulties.

A special case of Theorem 5.6 was given in \cite{Sk02}.

\medbreak
Although several fundamental facts of the classical theory generalize to 
commutative $H$-semiprime algebras, in the case when the Hopf algebra $H$ is 
not cocommutative it may not admit sufficiently many actions on commutative 
algebras. The next result has been obtained by Etingof and Walton \cite{Et-W14} 
when either $\chr\bbk=0$ or $\chr\bbk>0$ and $H$ is also semisimple. Its extension 
to the case when $H$ is not necessarily semisimple has been given in 
\cite{Sk-PAMS17}.

\proclaim
Theorem 5.7.
Assume\/ $\bbk$ to be algebraically closed. Then any action of a finite 
dimensional cosemisimple Hopf algebra $H$ on a commutative domain $A$ factors 
through an action of a group algebra, i.e. there exists a Hopf ideal $I$ of 
$H$ such that $I$ annihilates $A$ and $H/I$ is spanned by grouplike elements.  
\endproclaim

Etingof and Walton say that the action of $H$ on $A$ is \emph{inner faithful}\/ 
if $A$ is not annihilated by any nonzero Hopf ideal of $H$. In \cite{Et-W15}, 
\cite{Et-W16} they investigated the question of the existence of inner 
faithful actions on commutative domains for pointed Hopf algebras. Some 
pointed Hopf algebras admit such actions, while the others do not.

The fact that the annihilator of $A$ in $H$ is often nontrivial had been 
recognized much earlier. Cohen and Westreich pointed out in \cite{Coh-W93, 
Cor. 0.12} that $H$ can act faithfully (in the ordinary sense) on a field $A$ 
only if $H$ is involutory and all grouplikes of $H^*$ lie in the center of 
$H^*$. Here $A$ can be even a domain since then the action of $H$ extends to 
the quotient field.

All this shows that the class of commutative $H$-module algebras is too narrow 
when $H$ is not cocommutative, and there is a definite need to study the 
invariants in the larger class of algebras satisfying a polynomial identity. 
As yet, not all results known for commutative $H$-module algebras have been 
extended to the PI case however.

As an extension of the commutative theory in a different direction Cohen and 
Westreich \cite{Coh-W94} introduced \emph{quantum commutative} $H$-module 
algebras. The commutativity law in these algebras comes from the braiding 
determined by a quasitriangular structure on $H$. Cohen, Westreich and Zhu proved

\proclaim
Theorem 5.8 \cite{Coh-WZ96}.
Let $A$ be a quantum commutative $H$-module algebra where $H$ is triangular 
semisimple and either\/ $\chr\bbk=0$ or\/ $\chr\bbk>\dim H$. Then $A$ is integral 
over $A^H$ and $A$ is PI.
\endproclaim

One may wonder whether the conclusion of this theorem is valid under less 
stringent restrictions on $H$ and the characteristic of $\bbk$ when $A$ is 
$H$-semiprime.

\section
6. The $H$-equivariant Martindale quotient ring

Here we present results of Eryashkin \cite{Er-a16} on quotient rings of 
$H$-semiprime PI algebras. Generalized Martindale quotient rings can be 
defined with respect to any filter $\F$ of ideals of a ring $R$ subject to the 
conditions that each ideal $I\in\nobreak\F$ has zero left and right 
annihilators in $R$ and that $IJ\in\F$ whenever $I,J\in\F$. Details of this 
construction are given, e.g., in \cite{Mo93, \S6.4}. If $R$ is prime and $\F$ 
is the set of all nonzero ideals of $R$, this construction gives the 
\emph{left, right and symmetric Martindale rings of quotients}, as defined in 
\cite{Pas, Ch. 3}.

Let $A$ be an $H$-module algebra. Denote by $\FHA$ the set of all its $H$-stable 
ideals with zero left and right annihilators in $A$. If $A$ is $H$-prime, then 
$\FHA$ consists of all nonzero $H$-stable ideals of $A$. The Martindale quotient 
rings with respect to this filter were introduced by Cohen \cite{Coh86}. The 
use of $H$-stable ideals in their definition accounts for the extension of the 
$H$-module structure to these quotient rings.

We will be concerned only with the \emph{$H$-symmetric ring of quotients} 
$Q_H(A)$ (see \cite{Mo93, p. 98}). The larger left and right quotient rings 
are less useful. The ring $Q_H(A)$ is characterized by the following properties 
(cf. \cite{Pas, Prop. 10.4}):

\setitemsize(3)
\item(1)
$Q_H(A)$ contains $A$ as a subring;

\item(2)
each $I\in\FHA$ has zero left and right annihilators in $Q_H(A)$;

\item(3)
for each $q\in Q_H(A)$ there exists $I\in\FHA$ such that $Iq\sbs A$ and 
$qI\sbs A$;

\item(4)
given $I\in\FHA$, a left $A$-linear map $f_l:I\to A$, and a right $A$-linear 
map $f_r:I\to A$ such that $xf_r(y)=f_l(x)y$ for all $x,y\in I$, there exists  
$q\in Q_H(A)$ such that $f_l(x)=xq$ and $f_r(x)=qx$ for all $x\in I$.

\medskip
Put $Q=Q_H(A)$. As explained in \cite{Coh86}, $Q$ is an $H$-module algebra, and 
$A$ embeds in $Q$ as an $H$-stable subalgebra. The centers $Z(A)$, $Z(Q)$ of 
$A$ and $Q$ are not stable, in general, under the action of $H$. Set
$$
\ZAH=Z(A)\cap A^H,\qquad\ZQH=Z(Q)\cap Q^H.
$$
It follows from (2) and (3) in the characterization of $Q_H(A)$ above that each 
nonzero left or right $A$-submodule of $Q$ has nonzero intersection with $A$. 
Therefore $Q$ has to be $H$-prime or $H$-semiprime whenever so is $A$.

There are several general properties of the $H$-equivariant Martindale 
quotient rings of $H$-prime algebras. In particular, the fact that $\ZQH$ is a 
field was explicitly stated in Matczuk's paper \cite{Matc91, Lemma 1.4} which 
used the right quotient rings, however. This field is called the 
\emph{$H$-extended centroid}\/ of $A$.

\proclaim
Lemma 6.1.
Suppose that $A$ is $H$-prime. Then $\ZQH$ is a field. Furthermore, given any 
$H$-stable ideal $I$ of $Q$ and a morphism $f:I\to Q$ in $\HQMQ$, there exists 
$z\in\ZQH$ such that $f(x)=zx$ for all $x\in I$.
\endproclaim

\Proof.
Recall that $f$ is a homomorphism of $Q$-bimodules and $H$-modules. 
Put $I'=f^{-1}(A)\cap A$. This is an $H$-stable ideal of $A$. We may assume 
that $f\ne0$. Then $I'\ne0$, i.e. $I'\in\FHA$. Note that $f(x)y=f(xy)=xf(y)$ 
for all $x,y\in I$. By (4) there exists $z\in Q$ such that $f(x)=zx=xz$ for 
all $x\in I'$.

If $u\in Q$ is any element, then $uJ\sbs A$ for some $J\in\FHA$. Replacing $J$ 
with $JI'$, we will also have $J\sbs I'$ and $uJ\sbs I'$. Since $z$ centralizers 
all elements of $I'$, it follows that $(zu-uz)J=0$. Hence $zu=uz$ by (2). 

Since $f$ is $H$-linear, we deduce that
$$
(hz)x=\sum h\1f\bigl(S(h\2)x\bigr)=\ep(h)f(x)=\ep(h)zx\qquad
\hbox{for $h\in H$, $x\in I'$}.
$$
In other words, $\bigl(hz-\ep(h)z\bigr)I'=0$. Hence $hz=\ep(h)z$ for all $h\in 
H$, and we conclude that $z\in\ZQH$.

Define $g:I\to Q$ by the formula $g(x)=f(x)-zx$. Then $g|_{I'}=0$. If $u\in I$ 
and $J\in\FHA$ is such that $J\sbs I'$ and $uJ\sbs I'$, then $g(u)J=g(uJ)=0$, 
yielding $g(u)=0$. Hence $g=0$, and so $f(x)=zx$ for all $x\in I$.

The annihilator of $z$ in $Q$ is an $H$-stable ideal. If this ideal is 
nonzero, then property (2) entails $z=0$. Hence $f$ has to be injective 
whenever $z\ne0$. In this case $f(I)$ is a nonzero $H$-stable ideal of $Q$. 
Applying the already proved conclusion to the inverse map $f^{-1}:f(I)\to I$, 
we see that there exists $z'\in\ZQH$ such that $f^{-1}(y)=z'y$ for all 
$y\in f(I)$. Then $(zz'-1)I=0$, and it follows that $z'=z^{-1}$. But this 
argument applies to each nonzero element of the commutative ring $\ZQH$. 
Thus $\ZQH$ is a field.
\endproof

\proclaim
Lemma 6.2.
Suppose that $A$ is $H$-prime. Let $S$ be any simple algebra whose center 
contains $\ZQH$. Consider $S\ot_\ZQH Q$ as an $H$-module algebra with respect 
to the action of $H$ on the second tensorand. If $I$ is a nonzero $H$-stable 
ideal of this algebra, then $I$ has nonzero intersection with the image of $Q$ 
in $S\ot_\ZQH Q$ under the map $x\mapsto1\ot x$.
\endproclaim

\Proof.
Let $n$ be minimal possible for which $I$ contains an element $u\ne0$ which can 
be written as $u=a_1\ot b_1+\ldots+a_n\ot b_n$ with $a_i\in S$, $b_i\in Q$. Fix 
such an element and its expression as a sum. Put
$$
M=\{(x_1,\ldots,x_n)\in Q^n\mid\tsum a_i\ot x_i\in I\}.
$$
Consider $Q^n$ as an object of $\HQMQ$ with respect to the natural actions of 
$H$ and $Q$ on each component. Then $M$ is a subobject of $Q^n$ in this 
category. Note that
$$
(b_1,\ldots,b_n)\in M.
$$
Let $p_i:Q^n\to Q$, $\,i=1,\ldots,n$, be the projections. Then $J=p_1(M)$
is an $H$-stable ideal of $Q$. But $\Ker p_1|_M=0$ since otherwise $I$ would 
contain a nonzero element written as $\sum_{i=2}^na_i\ot x_i$ with less than 
$n$ summands. Thus $p_1|_M:M\to J$ is an isomorphism in the category 
$\HQMQ$. Setting $f_i=p_i\circ(p_1|_M)^{-1}$, we get
$$
M=\{\bigl(f_1(x),\ldots,f_n(x)\bigr)\mid x\in J\},
$$
and each map $f_i:J\to Q$ is a morphism in $\HQMQ$. By Lemma 6.1 there 
exist $z_1,\ldots,z_n\in\ZQH$ such that $f_i(x)=z_ix$ for all $x\in J$. In 
particular, $b_i=z_ib_1$ for each $i$, and therefore
$u=\sum a_i\ot z_ib_1=(\sum a_iz_i)\ot b_1$.

The minimality of $n$ implies that $n=1$, i.e. $u=a_1\ot b_1$. But then 
$Sa_1S\ot b_1\sbs I$. Since $S$ is simple, we have $Sa_1S=S$. Hence 
$1\ot b_1\in I$.
\endproof

\proclaim
Lemma 6.3.
Suppose that $A$ is $H$-prime, and let $P$ be any prime ideal of $A$ such that 
$P_H=0$. Denote by $\Qov$ the symmetric Martindale quotient ring of the prime 
ring $\Aov=A/P$. The canonical map $\pi:A\to\Aov$ extends to a ring homomorphism 
$Q\to\Qov$ which maps the center of $Q$ into the center of $\Qov$.
\endproclaim

\Proof.
Let $q\in Q$. There exists $I\in\FHA$ such that $Iq\sbs A$ and $qI\sbs A$. 
Since $P_H=0$, we have $I\not\sbs P$. Hence $\pi(I)$ is a nonzero ideal of 
$\Aov$. Note that $Iq(I\cap P)$ is contained in $P$. Applying $\pi$, we get 
$\pi(I)\pi\bigl(q(I\cap P)\bigr)=0$, whence $\pi\bigl(q(I\cap P)\bigr)=0$ 
since $\Aov$ is prime.

This shows that $q(I\cap P)\sbs P$. Similarly $(I\cap P)q\sbs P$. Therefore 
the right and left multiplications by $q$ induce, respectively, a left 
$\Aov$-linear map $f_l:\pi(I)\to\Aov$ and a right $\Aov$-linear map 
$f_r:\pi(I)\to\Aov$. The pair $(f_l,f_r)$ determines an element 
$\qov\in\Qov$. It is easy to see that the assignment $q\mapsto\qov$ defines 
a ring homomorphism $Q\to\Qov$ whose restriction to $A$ is $\pi$.

Denote this extension of $\pi$ by the same letter $\pi$. If $z\in Z(Q)$, then 
$\pi(z)$ commutes with all elements of $\pi(A)=\Aov$, but then $\pi(z)$ 
commutes with all elements of $\Qov$.
\endproof

\proclaim
Theorem 6.4 \cite{Er-a16}.
Suppose that $A$ is PI and $H$-prime. Then the $H$-symmetric quotient ring 
$Q=Q_H(A)$ is an $H$-simple $H$-module algebra of finite dimension over $\ZQH$. 
Moreover, $Q=\ZQH A$.
\endproclaim

\Proof.
Take any prime ideal $P$ of $A$ such that $P_H=0$. Let $\pi:Q\to\Qov$ be the 
ring homomorphism of Lemma 6.3. Since $\Aov$ is a prime PI algebra, $\Qov$ is 
the classical quotient ring of $\Aov$ (see \cite{Pas, Th. 23.4}). By Posner's 
theorem the ring $\Qov$ is simple and finite dimensional over its center 
$Z(\Qov)$. The composite map
$$
\ph:Q\lmapr2{\rho}Q\ot H^*\lmapr4{\pi\ot\id}\Qov\ot H^*
$$
is a homomorphism of $H$-module algebras, assuming the trivial action of 
$H$ on $\Qov$ and the hit action $\rhu$ on $H^*$. It extends to a homomorphism 
of $H$-module algebras
$$
\psi:Z(\Qov)\ot_{\ZQH}Q\to\Qov\ot H^*,\qquad z\ot q\mapsto(z\ot1)\cdot\ph(q).
$$
Since $(\id\ot\ep)\circ\ph=\pi$, we have $\Ker\ph\sbs\Ker\pi$. It follows that 
$$ A\cap\Ker\ph\sbs A\cap\Ker\pi=\Ker\pi|_A=P.
$$
But $\Ker\ph$ is an $H$-stable ideal of $Q$. Since $P_H=0$, we get 
$A\cap\Ker\ph=0$, which entails $\Ker\ph=0$. Now $\Ker\psi$ is an $H$-stable 
ideal of $Z(\Qov)\ot_{\ZQH}Q$. It has zero intersection with the image of $Q$ 
by the preceding conclusion, whence $\Ker\psi=0$ by Lemma 6.2. 
Injectivity of $\psi$ entails an upper bound for the dimension
$$
\eqalign{
[Q:\ZQH]&{}=[Z(\Qov)\ot_{\ZQH}Q:Z(\Qov)]\cr
&{}\le[\Qov\ot H^*:Z(\Qov)]=[\Qov:Z(\Qov)]\cdot\dim_\bbk H<\infty.\cr
}
$$
The $H$-stable subalgebra $A'=\ZQH A\sbs Q$ is then finite dimensional over 
$\ZQH$ too. Also, $A'$ is $H$-prime since each nonzero ideal of $A'$ has 
nonzero intersection with $A$. By Theorem 2.4 $A'$ is $H$-simple. But for each 
$q\in Q$ there exists a nonzero $H$-stable ideal $I'$ of $A'$ such that 
$qI'\sbs A'$. We must have $1\in I'$, and so $q\in A'$. Thus $Q=A'$.
\endproof

\proclaim
Corollary 6.5.
If $A$ is PI and $H$-prime, then $A$ has finitely many minimal prime ideals, 
and $P_H=0$ for each of them.
\endproclaim

\Proof.
Since $Q$ is artinian, it has finitely many maximal ideals. Let $P_1,\ldots,P_n$ 
be their contractions to $A$. The intersection $\bigcap P_i$ is nilpotent since 
it is contained in the Jacobson radical of $Q$. Hence each prime ideal of $A$ 
contains $P_i$ for some $i$, i.e. all minimal primes are among $P_1,\ldots,P_n$. 
If $I$ is any $H$-stable ideal of $A$ contained in $P_i$, then $IQ$ is an 
$H$-stable ideal of $Q$ contained in a maximal ideal. It follows that $IQ=0$, 
and therefore $I=0$.
\endproof

\proclaim
Corollary 6.6.
Suppose that $A$ is PI and $H$-prime. If $P$ is any prime ideal of $A$ such 
that $P_H=0$, then the ring homomorphism $\pi:Q\to\Qov$ of Lemma\/ {\rm6.3} 
is surjective.
\endproclaim

\Proof.
We have $\Aov=\pi(A)\sbs\pi(Q)\sbs\Qov$.
If $s$ is any regular element of $\Aov$, then $s$ is invertible in $\Qov$ 
since $\Qov$ is a classical quotient ring of $\Aov$. But $s\in\pi(Q)$, whence 
$s$ is a regular element of $\pi(Q)$. Since $\pi(Q)$ is a finite dimensional 
algebra over a field, it follows that $s^{-1}\in\pi(Q)$. Then $\Qov=\pi(Q)$.  
\endproof

\proclaim
Corollary 6.7.
If $A$ is PI and $H$-simple, then $A$ has finite dimension over its central 
subfield $Z(A)^H$.
\endproclaim

\Proof.
In this case $Q_H(A)=A$.
\endproof

\proclaim
Lemma 6.8.
Suppose that $K_1,\ldots,K_n$ are minimal $H$-prime ideals of $A$ such that 
$\bigcap K_i=0$. Then $\,Q_H(A)\cong Q_H(A/K_1)\times\ldots\times Q_H(A/K_n)$.  
\endproclaim

\Proof.
Put $Q=Q_H(A)$, $\,A_i=A/K_i$ and $Q_i=Q_H(A_i)$ for each $i$. The canonical map 
$\pi_i:A\to A_i$ extends to a homomorphism of $H$-module algebras $Q\to Q_i$ 
as in the proof of Lemma 6.3. The main point here is that $A_i$ is 
$H$-prime and $I\not\sbs K_i$ for each $I\in\FHA$. In fact $K_i$ has nonzero 
annihilator in $A$ since $K_iK'_i\sbs K_i\cap K'_i=0$ where 
$K'_i=\bigcap_{j\ne i}K_j$. Hence each element of $Q$ gives rise to a left 
$A_i$-linear map $f_l:\pi_i(I)\to A_i$ and a right $A_i$-linear map 
$f_r:\pi_i(I)\to A_i$, and the pair $(f_l,f_r)$ determines an element of $Q_i$.

Now the collection $\pi_1,\ldots,\pi_n$ gives a homomorphism of $H$-module 
algebras
$$
\pi:Q\to Q_1\times\ldots\times Q_n.
$$
Since $\Ker\pi|_A=\bigcap K_i=0$, it follows that $\Ker\pi=0$. It remains to 
show that $\pi$ is surjective. Suppose $q_1\in Q_1$. There exists a nonzero 
$H$-stable ideal $I_1$ of $A_1$ such that $I_1q_1\sbs A_1$ and $q_1I_1\sbs A_1$. 
Take any $H$-stable ideal $J$ of $A$ with the property that 
$0\ne\pi_1(J)\sbs I_1$. Replacing $J$ with $JK'_1$ we will also have
$$
J\sbs K'_1,\qquad\pi_1(J)\,q_1\sbs\pi_1(K'_1),\qquad 
q_1\,\pi_1(J)\sbs\pi_1(K'_1).
$$
Since $K_1\cap K'_1=0$, there is an isomorphism of $A$-bimodules 
$K'_1\cong\pi_1(K'_1)$. Define maps $\,f_l,f_r:J\to K'_1\,$ by the rules
$$
\pi_1\bigl(f_l(x)\bigr)=\pi_1(x)q_1,\qquad
\pi_1\bigl(f_r(x)\bigr)=q_1\pi_1(x).
$$
Put $I=J+\sum_{i\ne1}K'_i$. Note that $K'_i\sbs K_j$ for $j\ne i$, while 
$K'_j\cap K_j=0$. Hence the sum $\sum K'_i$ is direct, and there are 
extensions of $f_l,f_r$ to maps $I\to A$ vanishing on $K'_i$ for each $i\ne1$. 
Note that $f_l$ is left $A$-linear, while $f_r$ is right $A$-linear. 

Since $\pi_i(I)\ne0$ for each $i$, the left and right annihilators of $I$ in 
$A$ are contained in each $K_i$. Then these annihilators are zero, i.e. 
$I\in\FHA$. Hence the pair $(f_l,f_r)$ determines an element $q\in Q$ such that 
$\pi_1(q)=q_1$ and $\pi_i(q)=0$ for $i\ne1$. By symmetry $Q_1$ in this argument 
can be replaced with $Q_j$ for any $j$.
\endproof

\proclaim
Corollary 6.9.
Suppose that $A$ is PI and $H$-semiprime with finitely many minimal $H$-prime 
ideals. Then $Q_H(A)\cong Q_1\times\ldots\times Q_n$ where $Q_1,\ldots,Q_n$ are 
$H$-simple $H$-module algebras with $[Q_i:\ZQiH]<\infty$ for each $i$.  
\endproclaim

\Proof.
In any $H$-semiprime algebra the intersection of all minimal $H$-prime ideals 
is zero. So Lemma 6.8 applies. It gives a direct product decomposition of 
$Q_H(A)$ in which each factor $Q_i=Q_H(A/K_i)$ is $H$-simple and finite 
dimensional over $\ZQiH$ by Theorem 6.4.
\endproof

\proclaim
Theorem 6.10.
Suppose that $A$ is PI and $H$-semiprime with finitely many minimal $H$-prime 
ideals. Then $Q_H(A)$ is a classical right quotient ring of $A$.
\endproclaim

We have made some comments about the proof in the discussion following the 
statement of Theorem 2.12. Eryashkin proves the conclusion of Theorem 6.10 in 
the $H$-prime case. If $A$ is PI and $H$-semiprime, and if $K_1,\ldots,K_n$ 
are all its minimal $H$-prime ideals, then, knowing that $Q_i=Q_H(A/K_i)$ is a 
classical quotient ring of $A/K_i$ for each $i$, he concludes that 
$Q_1\times\ldots\times Q_n$ is a classical quotient ring of $A$ directly, not 
using Lemma 6.8.

For a special class of $H$-prime PI algebras Theorems 6.4 and 6.10 were 
obtained in an earlier article \cite{Er12}.

\proclaim
Lemma 6.11.
The set of minimal $H$-prime ideals of $A$ is finite if $A$ is finitely 
generated and PI. The same conclusion holds if $A$ is either left or right 
noetherian.
\endproclaim

\Proof.
Under each of these assumptions $A$ has finitely many minimal prime ideals and 
the prime radical $N$ of $A$ is nilpotent (in the case of a finitely generated 
PI algebra see \cite{Row, Cor. 6.3.36$'$, Th. 6.3.39}). Let $P_1,\ldots,P_n$ 
be all the minimal primes, and for each $i$ let $K_i$ be the largest 
$H$-stable ideal of $A$ contained in $P_i$. Since $\prod K_i\sbs N$ is 
nilpotent, any minimal $H$-prime ideal of $A$ has to coincide with one of the 
$H$-prime ideals $K_1,\ldots,K_n$.
\endproof

\section
7. Integrality of PI algebras over the invariants

If $A$ is a noncommutative $H$-module algebra, it is meaningful to consider 
integrality of $A$ over invariants in two different senses. One question 
concerns integrality over central invariants. For this it should be assumed 
at least that $A$ is integral over its center $Z(A)$. We denote by $\ZAH$ the 
subalgebra of $Z(A)$ consisting of $H$-invariant central elements.

If $H$ is cocommutative, then $Z(A)$ is stable under the action of $H$, and 
$Z(A)$ is integral over $\ZAH$ by the classical theory. If $H$ is not 
cocommutative, the problem becomes highly nontrivial. When $H$ is pointed, or 
at least when the coradical of $H$ is cocommutative the following result was 
obtained by Totok:

\proclaim
Theorem 7.1 \cite{Tot98}.
The center $Z(A)$ is integral over $\ZAH$, and therefore $A$ is integral over 
$\ZAH$ if $A$ is integral over $Z(A)$, in each of the following two cases\/{\rm:}

\item(a)
\quad$\chr\bbk>0$ and $H$ has cocommutative coradical,

\item(b)
\quad$\chr\bbk=0$, $H$ is pointed, $A$ is reduced, $Z(A)$ is finitely generated.

\endproclaim

Making use of the coradical filtration $H_0\sbs H_1\sbs\cdots\sbs H_n=H$, 
Totok constructs a chain of subalgebras $Z_0\sps Z_1\sps\cdots\sps Z_n$ such 
that $Z_0=Z(A)^{H_0}$, and for each $i>0$ the ring $Z_{i-1}$ is integral over 
$Z_i$ and $H_i$ acts trivially on $Z_i$ in the sense that $hz=\ep(h)z$ for all 
$h\in H_i$ and $z\in Z_i$. The conclusion of Theorem 7.1 follows then by 
transitivity of integrality since $Z(A)$ is integral over $Z_0$. This extends 
the technique applied by Artamonov \cite{Art96} in the case when $A$ is a 
commutative domain.

New results on integrality of $A$ over $\ZAH$ for an arbitrary finite 
dimensional Hopf algebra have appeared very recently. In \cite{Et15} Etingof 
observes that the problem admits a bimodule reformulation which can be studied 
independently of any Hopf algebra theory. In fact $\ZAH$ consists precisely 
of those elements $a\in A$ for which the left multiplication by $a\ot1$ in the 
algebra $A\ot H^*$ coincides with the right multiplication by $\rho(a)$ or, in 
other words, the left action of $a$ on $A\ot H^*$ is the same as the right 
action with respect to the $A$-bimodule structure defined as in section 3. This 
bimodule has a special property which has been used by Etingof to introduce 
the notion of \emph{Galois bimodules}.

An $R$-bimodule $P$ for a ring $R$ is called Galois of rank $d\ge1$ if $P$ is 
left and right free of rank $d$ and there is an isomorphism of bimodules 
$P\ot_RP\cong P^d$. Etingof derives a classification of Galois bimodules when 
$R$ is a semisimple artinian ring module-finite over its center $Z(R)$. Let 
$R\cong R_1\times\ldots\times R_n$ where $R_1,\ldots,R_n$ are simple rings. In 
the process of the classification it is verified that for each Galois bimodule 
$P$ the ring $R$ is a finite module over the \emph{center} of $P$ defined as
$$
Z(P)=\{a\in R\mid ax=xa\hbox{ for all $x\in P$}\}\sbs Z(R).
$$ 
Let $\phi_i(a)$ be the $R_i$-linear endomorphism of $R_i\ot_RP$ afforded by the 
right action of $a\in R$. Now $R_i\ot_RP$ is a finite dimensional vector space 
over the center $Z(R_i)$ of $R_i$, and $\phi_i(a)$ may be regarded as a linear 
transformation of this vector space. So the characteristic polynomial 
$\chi_{\phi_i(a)}\in Z(R_i)[t]$ makes sense. Let
$$
\chi_a\in Z(R)[t]\cong Z(R_1)[t]\times\ldots\times Z(R_n)[t]
$$
be the polynomial whose $i$th component is $\chi_{\phi_i(a)}^{m^2/m_i^2}$ 
where $m_i^2=[R_i:Z(R_i)]$ and $m$ is the least common multiple of 
$m_1,\ldots,m_n$. It is shown in \cite{Et15} that for each central element 
$a\in Z(R)$ all coefficients of $\chi_a$ belong to $Z(P)$. This is a key fact 
needed for applications to integrality.

Suppose that the $H$-module algebra $A$ has a semisimple artinian 
classical quotient ring $Q$ which is a finite module over its center $Z(Q)$. 
Then Etingof's results discussed in the preceding paragraph apply to the 
Galois $Q$-bimodule $Q\ot H^*$ whose center is $\ZQH$. In particular, $Q$ is 
module-finite over $\ZQH$ and certain polynomials associated with central 
elements of $Q$ have coefficients in $\ZQH$.

However, the ultimate goal is to find conditions ensuring integrality of $A$. 
Etingof formulates results for comodule algebras, but we stick to the 
conventions set for the present paper.

\proclaim
Theorem 7.2 \cite{Et15}.
Let $Z$ be a central subalgebra of $A$ whose total quotient ring $Q(Z)$ is a 
direct product of finitely many fields. Suppose that

\item(1)
$A$ is a finitely generated torsion-free $Z$-module,

\item(2)
$Q(Z)\otimes_ZA$ is a semisimple ring with center $Q(Z)$,

\item(3)
either $Z$ is integrally closed in $Q(Z)$ or $A$ is a projective $Z$-module.

\smallskip
\noindent
Then $A$ is integral over $Z\cap A^H$.
\endproclaim

An $H$-module algebra $A$ is said to be \emph{indecomposable} if it is not 
isomorphic to a direct product of two nonzero $H$-module algebras. If $A$ is 
artinian and $H$-semiprime, this is equivalent, by Theorem 2.4, to the 
$H$-simplicity of $A$.

Indecomposability of $A$ in the next proposition means that the corresponding 
Galois $A$-bimodule $P=A\ot H^*$ is connected. In this case $P^{m_*^2}$ is 
isomorphic to a multiple of $A\ot_{Z(P)}A$ by the classification of Galois 
bimodules. Here $Z(P)=\nobreak\ZAH\!$. Comparing the left ranks of the two 
bimodules Etingof deduces a divisibility relation involving numeric 
characteristics of $A$:

\proclaim
Proposition 7.3 \cite{Et15}.
Suppose that $A$ is semisimple artinian, module-finite over $Z(A)$, and 
indecomposable as an $H$-module algebra. Let $A_1,\ldots,A_n$ be all simple 
factor rings of $A$. Put
$$
d_i=[Z(A_i):\ZAH],\quad m_i=[A_i:Z(A_i)]^{\mskip1mu1/2},\quad
m_*=\gcd(m_1,\ldots,m_n).
$$
Then $\sum d_i(m_i/m_*)^2$ divides the dimension of $H$. In other words, 
$\,[A:\ZAH]$\hfil\break divides $\,m_*^2\,(\dim H)$.
\endproclaim

Now we describe a different approach, due to Eryashkin, which makes 
systematic use of structural properties of $H$-module algebras discussed 
earlier. The results have been obtained not only for semiprime $H$-module 
algebras but for $H$-semiprime algebras as well.

The ring homomorphism $A\to A\ot H^*$ given by the assignment $a\mapsto a\ot1$ 
maps the center of $A$ into a central subalgebra of $A\ot H^*$. Therefore 
$A\ot H^*$ may be regarded as a $Z$-algebra for any central subalgebra $Z$ of 
$A$. If $A$ is projective of finite constant rank as a $Z$-module, then so too 
is $A\ot H^*$. As explained in section 5, in this case there are characteristic 
polynomials for the ring extension $A\ot H^*/Z$. For each $a\in A$
$$
P_{A\ot H^*/Z}\bigl(\rho(a),t\bigr)\in Z[t]
$$
is the characteristic polynomial of the left multiplication operator by the 
element $\rho(a)$ in the $Z$-algebra $A\ot H^*$. Alternatively, one could use 
the characteristic polynomials of the right multiplication operators.

At one point we will need the ring-theoretic fact stated below. For the proof 
see \cite{Er-a16, Prop. 3.1}.

\proclaim
Proposition 7.4 \cite{Er-a16}.
Let $R$ be a ring which has a right artinian classical right quotient ring 
$Q(R)$. Suppose that $R$ is a finitely generated module over a central subring 
$Z$ such that $\,\ann_R(z)=\ann_Z(z)R\,$ for each $z\in Z$. Then $Q(R)\cong 
Q(Z)\ot_ZR$ where $Q(Z)$ is the total quotient ring of $Z$.
\endproclaim

\proclaim
Theorem 7.5 \cite{Er-a16}.
Suppose that $A$ is $H$-semiprime with finitely many minimal \hbox{$H$-prime} 
ideals. If $A$ is projective of finite constant rank as a module over its 
center $Z(A)$, then $A$ is integral over $\ZAH$. In fact, for each $a\in A$ 
the characteristic polynomial $\,P_{A\ot H^*\!/Z(A)}\bigl(\rho(a),t\bigr)\,$ 
has all coefficients in $\ZAH$.
\endproclaim

\Proof.
Before we treat the general case let us verify the conclusion of this theorem 
under additional assumptions about $A$.

\medbreak
Step 1. Suppose that $A$ is $H$-simple.

By Corollary 6.7 $\ZAH$ is a field and $[A:\ZAH]<\infty$. Integrality of $A$ 
over $\ZAH$ is immediate. We still have to prove the statement about the 
characteristic polynomials.

Consider the action of $H$ on $A\ot H^*$ defined by the rule 
$h(a\ot\xi)=a\ot(h\rhu\xi)$ for $h\in H$, $a\in A$ and $\xi\in H^*$. Then 
the map $\rho:A\to A\ot H^*$ is a homomorphism of $H$-module algebras. So too 
is its extension
$$
\psi:Z(A)\ot_{\ZAH}A\to A\ot H^*,\qquad z\ot a\mapsto(z\ot1)\cdot\rho(a),
$$
with the action of $H$ on the first algebra defined by the formula 
$h(z\ot a)=z\ot ha$ for $h\in H$, $z\in Z(A)$ and $a\in A$. Put
$$
\Ascr=Z(A)\ot_{\ZAH}A,\qquad\Bscr=A\ot H^*.
$$
We claim that $\Bscr$ is a free left $\Ascr$-module with respect to the 
action afforded by $\psi$. Since the ring $Z(A)$ is artinian and since $\psi$ 
is a homomorphism of $Z(A)$-algebras, both of which are free modules over 
$Z(A)$, it suffices to check that for each maximal ideal $\m$ of $Z(A)$ the 
$\Ascr/\m\Ascr$-module $\Bscr/\m\Bscr$ is free of rank $r$ where $r$ does not 
depend on $\m$. Now
$$
\Ascr/\m\Ascr=K(\m)\ot_{\ZAH}A,\qquad\Bscr/\m\Bscr=A/\m A\ot H^*
$$
with $K(\m)=Z(A)/\m$ being a field. Since $A$ is $H$-simple, it follows from 
Lemma 6.2 that $\Ascr/\m\Ascr$ is $H$-simple too. Hence $\Bscr/\m\Bscr$ is a 
free $\Ascr/\m\Ascr$-module by Theorem 2.9. The rank $r(\m)$ of this free 
module can be computed as
$$
r(\m)={[\Bscr/\m\Bscr:K(\m)]\over[\Ascr/\m\Ascr:K(\m)]}
={[A:Z(A)]\cdot(\dim H)\over[A:\ZAH]}
={\dim H\over[Z(A):\ZAH]}
$$
where $[A:Z(A)]$ is the rank of $A$ as a $Z(A)$-module. This shows that 
$r(\m)$ has the same value for all $\m$, as required.

Thus $\Bscr\cong\Ascr^r$ as an $\Ascr$-module. Since $\rho(a)=\psi(1\ot a)$, 
we deduce that
$$
P_{\Bscr/Z(A)}\bigl(\rho(a),t\bigr)
=P_{\Ascr/Z(A)}(1\ot a,t)^r=P_{A/\ZAH}(a,t)^r,
$$
which is a polynomial with coefficients in $\ZAH$. 

\medbreak
Step 2. Suppose that $A$ is artinian and $H$-semiprime.

By Theorem 2.4 $A=A_1\times\ldots\times A_n$ where each $A_i$ is 
an $H$-simple $H$-module algebra. Clearly
$$
Z(A)=Z(A_1)\times\ldots\times Z(A_n)\qquad{\rm and}\qquad
A^H=A_1^H\times\ldots\times A_n^H.
$$
Let $\pi_i:A\to A_i$ be the projection, $\ze_i:Z(A)\to Z(A_i)$ the 
restriction of $\pi_i$, and $\ze_i^t:Z(A)[t]\to Z(A_i)[t]$ the extension of 
$\ze_i$ to polynomial rings. Note that
$$
A_i\cong Z(A_i)\ot_{Z(A)}A,\quad\hbox{and therefore}\quad
A_i\ot H^*\cong Z(A_i)\ot_{Z(A)}(A\ot H^*).
$$
Since $(\pi_i\ot\id)\rho(a)=\rho(\pi_ia)$, we get
$$
\ze_i^t\,P_{A\ot H^*/Z(A)}\bigl(\rho(a),t\bigr)
=P_{A_i\ot H^*/Z(A_i)}\bigl(\rho(\pi_ia),t\bigr)\in Z(A_i)^H\,[t]
$$
by Step 1. It follows that all coefficients of the polynomial 
$P_{A\ot H^*/Z(A)}\bigl(\rho(a),t\bigr)$ are $H$-invariant since they have 
$H$-invariant images in each $A_i$.

\medbreak
It is easy now to complete the proof of Theorem 7.5 in full generality. By 
Theorem 2.12 $A$ has a right and left artinian classical right quotient ring 
$Q=Q(A)$ which is an $H$-semiprime $H$-module algebra since so is $A$. Note 
that $\,\ann_A(z)=\ann_{Z(A)}(z)A\,$ for each $z\in Z(A)$ since $A$ is a direct 
summand of a free $Z(A)$-module. By Proposition 7.4 $Q$ is a central localization 
of $A$. Then the total quotient ring of $Z(A)$ coincides with the center of 
$Q$, and so $Q\cong Z(Q)\ot_{Z(A)}A$. From the functorial properties of 
characteristic polynomials it follows that
$$
P_{A\ot H^*/Z(A)}\bigl(\rho(a),t\bigr)=P_{Q\ot H^*/Z(Q)}\bigl(\rho(a),t\bigr).
$$
All coefficients of this polynomial lie in $\ZQH$ by Step 2. Hence they 
actually lie in $\,Z(A)\cap\ZQH=\ZAH$.
\endproof

All ideas of this proof are taken from \cite{Er-a16}. We have used Theorem 2.9 
to make some arguments more transparent. Note that Step 1 in the proof yields 
also the following conclusion:

\proclaim
Corollary 7.6.
Suppose that $A$ is $H$-simple and $A$ is a free module of finite rank over 
its center $Z(A)$. Then the dimension $[Z(A):\ZAH]$ of $Z(A)$ over $\ZAH$ 
divides the dimension of $H$.
\endproclaim

It is not clear to what extent the conditions imposed on $A$ in Theorem 7.5 
are optimal. One concern arising here is the finiteness of the set of minimal 
$H$-primes. An easy extension of Theorem 7.5 is stated below:

\proclaim
Corollary 7.7.
Suppose that $A$ is projective of finite constant rank as a module over its 
center $Z(A)$ and there is a set $\F$ of $H$-semiprime ideals of $A$ such that

\setitemsize(4)
\item(1)
each ideal in $\F$ is an intersection of finitely many $H$-prime ideals,

\item(2)
each ideal $I\in\F$ is generated by $I\cap Z(A)$,

\item(3)
for each $I\in\F$ the image of $Z(A)$ in $A/I$ coincides with the center 
of $A/I$,

\item(4)
$\bigcap_{I\in\F}I=0$.

\noindent
Then for each $a\in A$ the characteristic polynomial 
$\,P_{A\ot H^*/Z(A)}\bigl(\rho(a),t\bigr)\,$ has all\hfil\break 
coefficients in $\ZAH$.
\endproclaim

\Proof.
For each $I\in\F$ we have $A/I\cong Z(A/I)\ot_{Z(A)}A$, and this algebra 
satisfies the hypothesis of Theorem 7.5. Hence all coefficients of $P_{A\ot 
H^*/Z(A)}\bigl(\rho(a),t\bigr)$ have $H$-invariant images in $A/I$. But (4) 
ensures that any element $c\in A$ is $H$-invariant whenever $c+I$ is 
$H$-invariant in $A/I$ for each $I\in\F$.
\endproof

The assumptions about $A$ in Corollary 7.7 are admittedly too restrictive. 
However, they are satisfied when $A$ is commutative and $H$-semiprime. Thus 
Theorem 5.2  is a special case of Corollary 7.7.

Without projectivity of $A$ over $Z(A)$ the characteristic polynomials for the 
ring extension $A\ot H^*/Z(A)$ are not defined. One can still exploit 
finiteness of $Q=Q(A)$ over $\ZQH$. Recall that, if $A$ is PI and $H$-prime, 
then $Q$ is $H$-simple and $\ZQH$ is a field. The next result is based on 
\cite{Er-a16, Prop. 3.3}, although it was not stated this way.

\proclaim
Proposition 7.8.
Suppose that $A$ is $H$-prime and module-finite over its center $Z(A)$. Let 
$Q$ be the classical quotient ring of $A$.
For each $a\in A$ all coefficients of the characteristic polynomial 
$\,P_{Q/\ZQH}(a,t)\,$ are integral over $Z(A)$.
\endproclaim

\Proof.
Since $[Q:\ZQH]<\infty$, the set $\Max Q$ of all maximal ideals of $Q$ is 
finite. For each $M\in\Max Q$ denote by $\pi_M$ the canonical map $Q\to Q/M$. 
Here $Q/M$ is a simple ring finite dimensional over its center $Z(Q/M)$. Since 
$Z(A)\sbs Z(Q)$, we have $\pi_M\bigl(Z(A)\bigr)\sbs Z(Q/M)$.

An arbitrary element $q\in Q$ is integral over $Z(A)$ if and only if $\pi_M(q)$ 
is integral over $\pi_M\bigl(Z(A)\bigr)$ for each $M\in\Max Q$. Indeed, if the 
latter property holds, then for each $M$ there exists a polynomial $f_M$ in one 
indeterminate with all coefficients in $Z(A)$ and the leading coefficient 1 such 
that $f_M(q)\in M$. Putting 
$$
f=\prod_{M\in\Max Q}f_M,
$$
we will have $f(q)\in J$ where $J$ stands for the Jacobson radical of $Q$. But 
$J$ is nilpotent, whence $f(q)^n=0$ for some integer $n>0$. Clearly $f^n$ is a 
polynomial with all coefficients in $Z(A)$ and the leading coefficient 1.

We will check that the necessary and sufficient condition of integrality from 
the preceding paragraph is satisfied for all coefficients of 
$P_{Q/\ZQH}(a,t)$. This will yield the final conclusion.

Consider the composite map 
$\rho_M:Q\lmapr2\rho Q\ot H^*\lmapr5{\pi_M\ot\id}Q/M\ot H^*$ 
and its extension
$$
\psi:Z(Q/M)\ot_\ZQH Q\to Q/M\ot H^*,\qquad z\ot q\mapsto (z\ot1)\cdot\rho_M(q).
$$
With the $H$-module structures as in the proof of Theorem 7.5 $\psi$ is a 
homomorphism of $H$-module algebras, and the first algebra is $H$-simple by 
Lemma 6.2. It follows, in particular, that $\psi$ is injective. By Theorem 2.9 
$\psi$ makes $Q/M\ot H^*$ a free left module over $Z(Q/M)\ot_\ZQH Q$. Let $r$ 
be its rank (it depends on $M$).

Denote by $P_2$ and $P_1$ the characteristic 
polynomials of these two rings regarded as finite dimensional algebras over the 
field $Z(Q/M)$. Since $\rho_M(a)=\psi(1\ot a)$, we get
$$
P_2\bigl(\rho_M(a),t\bigr)=P_1(1\ot a,t)^r=P_{Q/\ZQH}(a,t)^r.
$$
Here we identify the field $\ZQH$ with its image in $Z(Q/M)$ under the map 
$\pi_M$.

Recall that $P_2\bigl(\rho_M(a),t\bigr)$ is the characteristic polynomial of 
the left multiplication operator associated with $\rho_M(a)$. But $\rho(a)\in 
A\ot H^*$, and therefore $\rho_M(a)$ lies in the subring $\pi_M(A)\ot H^*$ of 
$Q/M\ot H^*$ which is a finitely generated module over $\pi_M\bigl(Z(A)\bigr)$. 
It follows that $\rho_M(a)$ satisfies a polynomial relation of integral 
dependence over $\pi_M\bigl(Z(A)\bigr)$, whence so too does the corresponding 
multiplication operator. This means that all eigenvalues of this operator in 
any algebraic closure of $Z(Q/M)$ are integral over $\pi_M\bigl(Z(A)\bigr)$. 
But these eigenvalues are precisely the roots of the characteristic 
polynomial, i.e. the roots of $P_{Q/\ZQH}(a,t)$ in view of the equality above. 
The coefficients of $P_{Q/\ZQH}(a,t)$ are evaluations of the elementary 
symmetric functions at the roots of this polynomial. Hence they are integral 
over $\pi_M\bigl(Z(A)\bigr)$ too.
\endproof

\proclaim
Theorem 7.9 \cite{Er-a16}.
Suppose that $A$ is $H$-semiprime with finitely many minimal \hbox{$H$-prime} 
ideals. Let $Q$ be the classical quotient ring of $A$. If $Z(A)$ is integrally 
closed in $Z(Q)$ and $A$ is a finitely generated $Z(A)$-module, then $A$ is 
integral over $\ZAH$.
\endproclaim

\Proof.
By Theorem 2.4 $Q\cong Q_1\times\ldots\times Q_n$ where each $Q_i$ is an 
$H$-simple $H$-module algebra. Let $e_i\in Q$ be the $H$-invariant central 
idempotent whose projection to $Q_j$ is 1 for $j=i$ and 0 otherwise. Since 
$e_i\in Z(Q)$ is integral over $Z(A)$, we must have $e_i\in Z(A)$. Then 
$$
A\cong A_1\times\ldots\times A_n\quad\hbox{where $A_i=A/(1-e_i)A$}.
$$
Clearly $Q_i$ is the classical quotient ring of $A_i$. It follows that $A_i$ 
is $H$-prime, $A_i$ is module-finite over its center $Z(A_i)$, and $Z(A_i)$ is 
integrally closed in $Z(Q_i)$.

This reduces the proof to the case when $A$ is $H$-prime. But in this case 
Proposition 7.8 applies. It shows that for each $a\in A$ all coefficients of 
the characteristic polynomial $P_{Q/\ZQH}(a,t)$ are in $Z(A)\cap\ZQH=\ZAH$.
Since $P_{Q/\ZQH}(a,a)=0$ by a general property of characteristic polynomials, 
$a$ is integral over $\ZAH$.
\endproof

The statement of \cite{Er-a16, Prop. 3.3} contains the additional assumption 
that $\ann_A(z)$ is equal to $\ann_{Z(A)}(z)A$ for each $z\in Z(A)$. The proofs 
given above show that this assumption is not needed.

In connection with Theorems 7.2, 7.5, 7.9 we are prompted to ask

\proclaim
Question 7.10.
Is there any example of an $H$-semiprime algebra $A$ module-finite 
over its center $Z(A)$ such that $A$ is not integral over $\ZAH$\/{\rm?}
\endproclaim

\proclaim
Proposition 7.11 \cite{Er-b16}.
Suppose that $H$ is semisimple. If $A$ is PI and $A^H\sbs Z(A)$, then
$A$ is integral over $A^H$.
\endproclaim

\Proof.
It suffices to consider the case when $A$ is finitely generated. Then 
integrality of $A$ over $A^H$ is equivalent to module-finiteness. There exists 
a surjective homomorphism of $H$-module algebras $B\to A$ where $B$ has a 
grading $B=B_0\oplus B_1\oplus\ldots$ with finite dimensional $H$-stable 
homogeneous components such that $B_0=\bbk$ and $B_1$ generates $B$.
For this we can start with the tensor algebra of any finite dimensional 
$H$-submodule $V\sbs A$ which generates $A$ as an algebra. Taking the factor 
algebra of $B$ by a suitable $H$-stable ideal, we may assume that $B$ is PI. 
Factoring out another ideal generated by all commutators $xy-yx$ with $x\in B$ 
and $y\in B^H$, we may also assume that $B^H\sbs Z(B)$.

Since $H$ is semisimple, $B^H$ is mapped onto $A^H$. Therefore it suffices to 
show that $B$ is a finite $B^H$-module. By the graded Nakayama lemma this 
holds if and only if $\dim B/B^H_+B<\infty$ where $B^H_+=\sum_{i>0}B_i^H$.

\smallskip
Put $D=B/B^H_+B$. Note that $B^H_+B$ is a homogeneous $H$-stable ideal of 
$B$. Hence $D$ inherits the structure of a graded $H$-module algebra, and $D$ 
is PI. Since $B^H$ maps onto $D^H$, we have $D^H=\bbk$. It follows that 
$D_+^H=0$ where $D_+=\sum_{i>0}D_i$.

Let $P$ be any maximal ideal of $D$, and let $P_H$ be the largest $H$-stable 
ideal contained in $P$. By Kaplansky's theorem the simple algebra $D/P$ is 
finite dimensional over its center, and, since $D$ is finitely generated, we 
have $\dim D/P<\infty$ (over $\bbk$). Now $P_H$ is the kernel of the composite 
map $D\lmapr2\rho D\ot H^*\lmapr1{}D/P\ot H^*$. It follows that $\dim 
D/P_H<\infty$ too. By Theorem 2.4 $D/P_H$ is $H$-simple, which means that 
$P_H$ is a maximal $H$-stable ideal of $D$.

If $P_H\not\sbs D_+$, then $P_H+D_+=D$. Since all $H$-modules are completely 
reducible, we deduce that $P^H+D_+^H=D^H$ where $P^H=P\cap D^H$. Hence there 
exists $d\in D_+^H$ such that $d\notin P$. This entails $D_+^H\ne0$, a 
contradiction.

Thus $P_H\sbs D_+$ is the only possibility. Since $D_+$ is an $H$-stable 
ideal of $D$, we get $P_H=D_+$ by maximality of $P_H$, but then $P=D_+$ too. 
We conclude that $D$ has a single maximal ideal. Recall that the prime radical 
of any finitely generated PI algebra is nilpotent by the Braun theorem 
\cite{Row, Th. 6.3.39} and coincides with the intersection of all maximal 
ideals by the Amitsur-Procesi theorem \cite{Row, Th. 6.3.3}. This implies 
that $D_+$ is the prime radical of $D$ and that $D_+$ is nilpotent. But then 
$D_i=0$ for sufficiently large $i$. Hence $\dim D<\infty$, as required. 
\endproof

The condition $A^H\sbs Z(A)$ may look artificial, but sometimes it arises very 
naturally. For instance, this inclusion always holds when $A$ is quantum 
commutative. In \cite{Coh-W93} Cohen and Westreich 
investigated how the condition $A^H\sbs Z(A)$ affects various 
properties of an $H$-module algebra, especially in the case when $A$ is an 
$H^*$-Galois extension of $A^H$.

In \cite{Er11} and \cite{Er12} Eryashkin considered a special class $\Ascr$ of 
$H$-module algebras. An $H$-module algebra $A$ belongs to $\Ascr$ if $A$ 
has an ideal $I$ such that the factor algebra $A/I$ is commutative and $I$ 
contains no nonzero $H$-stable ideals of $A$. Such an algebra is PI since it 
embeds into the algebra $A/I\ot H^*$ which is a finite module over its center. 
If $z\in A^H$, then $\{za-az\mid a\in A\}$ is an $H$-submodule of $A$ contained 
in the ideal $I$, whence $za=az$ for all $a\in A$ by the conditions imposed on 
$I$. This shows that $A^H\sbs Z(A)$.

Starting with an arbitrary left $H$-module $V$ one obtains an $H$-prime 
algebra in $\Ascr$ taking $A=T(V)/I'_H$ where $T(V)$ is the tensor 
algebra of $V$ and $I'_H$ is its largest $H$-stable ideal contained in the 
ideal $I'$ generated by all commutators. Here the ideal $I=I'/I'_H$ of $A$ has 
the property that $A/I$ is the symmetric algebra of $V$. By a careful 
examination Eryashkin has verified that $A$ is not integral over $A^H$ in the 
case when $\chr\bbk=0$, $H$ is the 4-dimensional Hopf algebra described by 
Sweedler, and $V$ is one of its 2-dimensional indecomposable modules.

In particular, the semisimplicity of $H$ is necessary in Proposition 7.11, even 
if $A$ is assumed to be $H$-prime. In positive characteristic the previous 
construction does not give such an example (see Corollary 8.4). This leaves open

\proclaim
Question 7.12.
Suppose that $\chr\bbk>0$. Is there any example of an $H$-prime PI algebra 
$A$ such that $A^H\sbs Z(A)$, but $A$ is not integral over $A^H${\rm?}
\endproclaim

\medbreak
If $A^H\not\sbs Z(A)$, then integrality of $A$ over $A^H$ should be understood 
as defined by Schelter \cite{Sch76}. An element $x\in A$ is called 
\emph{Schelter integral} over $A^H$ if there exists an integer $m>0$ such that 
$x^m$ can be written as a sum of several elements, each of which is a product 
of elements contained in $A^H\cup\{x\}$ with $x$ occurring as a factor in this 
product less than $m$ times. If all elements of $A$ are Schelter integral over 
$A^H$, then $A$ is said to be Schelter integral over $A^H$.

In the 1993 expository lectures Montgomery asked whether Schelter integrality 
of $A$ over $A^H$ holds whenever $H$ is semisimple \cite{Mo93, Question 4.3.1}. 
By that time the group action case had been settled in full generality by 
Quinn. If $G$ is a finite automorphism group of a ring $R$ such that $|G|R=R$, 
then $R$ is Schelter integral over the subring of invariants $R^G$. In fact it 
was proved in \cite{Qu89} that $R$ is \emph{fully integral}\/ over $R^G$, which 
is a stronger property defined in terms of collections of elements rather than 
single elements. Quinn also obtained a partial result for Hopf actions:

\proclaim
Theorem 7.13 \cite{Qu89}.
Suppose that $H$ is semisimple and the action of $H$ on $A$ is inner. Then 
each ideal $I$ of $A$ is fully integral, and therefore also Schelter integral, 
over $I^H$ of degree bounded by a function in the dimension of $H$.
\endproclaim

The condition that the action is \emph{inner}\/ means that there exists an 
invertible element $u\in A\ot H^*$ such that
$$
\rho(a)=u(a\ot1)u^{-1}\quad\hbox{for all $a\in A$}.
$$
In particular, the two nonunital subalgebras $I\ot1$ and $\rho(I)$ of $A\ot H^*$ 
are conjugate by an inner automorphism. It follows that $(I\ot H^*)\#H\cong 
I\ot\Endk H$ is fully integral over $\rho(I)\#H$ since $I\ot\Endk H$ is known 
to be fully integral over $I\ot\bbk$ by the Par\'e-Schelter theorem 
\cite{Par-S78}. Finally, Quinn deduces that $I$ is fully integral over $I^H$ 
using an idempotent $e\in(A\ot H^*)\#H$ such that
$$
e\bigl((I\ot H^*)\#H\bigr)e\cong I\quad\hbox{and}\quad
e\bigl(\rho(I)\#H\bigr)e\cong I^H
$$
as nonunital algebras. In the case of inner action all ideals of $A$ 
are $H$-stable.

At present it is not known how to extend the previous theorem to arbitrary 
module algebras for a semisimple Hopf algebra. Special cases of the problem 
were dealt with in \cite{Coh86}, \cite{Tot98}. Eryashkin has succeeded in 
answering Montgomery's question in the case of PI algebras:

\proclaim
Theorem 7.14 \cite{Er-b16}.
Suppose that $H$ is semisimple and cosemisimple. If $A$ is PI, then
$A$ is Schelter integral over $A^H$.
\endproclaim

\Proof.
The initial idea comes from a paper of Montgomery and Small \cite{Mo-Sm84} 
where a similar problem for group actions on PI rings was considered. By 
Zorn's lemma $A$ has an $H$-stable ideal $K$ maximal with respect to the 
property that all elements of $K$ are Schelter integral over $A^H$. It is easy 
to see that $K$ is $H$-semiprime.

Replacing $A$ with $A/K$, we may assume that $A$ is $H$-semiprime and $A$ has 
no nonzero $H$-stable ideals which are Schelter integral over $A^H$. We have 
to show that $A=0$. This step requires more effort as compared with the group 
action case.

Suppose that $A\ne0$. Since $H$ is cosemisimple, Linchenko and Montgomery tell 
us that the prime radical of $A$ is $H$-stable \cite{Lin-M07, Th. 3.5}; hence 
$A$ is semiprime. Then, by the general PI theory, $A$ has a nonzero ideal 
$I'$ such that for each $x\in I'$ the left ideal $Ax$ is contained in a 
finitely generated $Z(A)$-submodule of $A$. Put $I=HI'$. This is a nonzero 
$H$-stable ideal of $A$. We will show that all elements of $I$ are Schelter 
integral over $A^H$, but this contradicts the assumptions about $A$.

Denote by $C$ the centralizer of $A^H$ in $A$. This is an $H$-stable subalgebra 
of $A$ with $C^H\sbs Z(C)$ and $Z(A)\sbs C$. Let $x\in I$. From the construction 
of $I$ it follows easily that $Ax$ is contained in a finitely generated 
$C$-submodule, say $N$, of $A$. Suppose that $N$ is generated as a $C$-module 
by $e_1,\ldots,e_n$. There exists a finitely generated subalgebra $C_0\sbs C$ 
such that $x$ and all elements $e_ix$ belong to the $C_0$-submodule $N_0$ of 
$N$ generated by $e_1,\ldots,e_n$. Then $N_0x\sbs N_0$.

Without loss of generality we may assume $C_0$ to be $H$-stable. Since 
$C_0^H\sbs Z(C_0)$, Proposition 7.11 shows that $C_0$ is integral, and therefore 
module-finite, over $C_0^H$. Hence $N_0$ is a finitely generated module over 
$C_0^H$. Define $r_x\in\End_{C_0}N_0$ by the rule $r_x(y)=yx$ for $y\in N_0$. 
Since $C_0^H$ is a commutative ring, the endomorphism $r_x$ satisfies a relation
$$
r_x^m+c_1r_x^{m-1}+\ldots+c_m\,\id=0
$$
for some integer $m>0$ and elements $c_1,\ldots,c_m\in C_0^H$. Applying this 
operator to $x\in N_0$, we deduce that $x^{m+1}+c_1x^m+\ldots+c_mx=0$.
It follows that $x$ is Schelter integral over $C_0^H\sbs A^H$.
\endproof

\section
8. Comparison with the invariants of the coradical

We continue to assume that $A$ is an $H$-module algebra. If $H$ is pointed with 
the group $G$ of grouplike elements, it was observed by Artamonov \cite{Art96} 
that $A^H=A^G$ when $\chr\bbk=0$ and $A$ is a commutative domain. If $\chr\bbk=p>0$, 
the Hopf algebra is pointed, and $A$ is commutative, then it follows from the 
results of Totok \cite{Tot97} and Zhu \cite{Zhu96} that $z^{p^s}\in A^H$ for 
all $z\in A^G$ where $s$ is the length of the coradical filtration of $H$. 

Etingof and Walton \cite{Et-W15} proved the equality $\ZAH=Z(A)^{H_0}$ where 
$H_0$ is the coradical of $H$ in the case when $A$ is a prime Azumaya algebra 
and $\chr\bbk=0$. In this section we present stronger conclusions, due to 
Eryashkin.

\proclaim
Proposition 8.1 \cite{Er-b16}.
Let $H_0\sbs H$ be a Hopf subalgebra containing the coradical of $H$. Suppose 
that $A$ is PI and $H$-simple. Let $K$ be a maximal $H_0$-stable ideal of $A$ 
and $A_0=A/K$. Denote by $\nu$ the canonical map $A\to A_0$.

\setitemsize(ii)
\item(i)
If $\chr\bbk=0$, then $\ZAcH=\nu\bigl(\ZAH\bigr)$.

\item(ii)
If $\chr\bbk=p>0$, then there exists an integer $s\ge0$ such that 
$z^{p^s}\in\nu\bigl(\ZAH\bigr)$ for all $z\in\ZAcH$.

\endproclaim

\Proof.
By Corollary 6.7 $\ZAH$ is a field and $[A:\ZAH]<\infty$. Similarly, $\ZAcH$ 
is a field and $[A_0:\ZAcH]<\infty$ since $A_0$ is PI and $H_0$-simple. Note 
that $\nu$ maps $\ZAH$ into $\ZAcH$. Define $H$-module structures on
$$ 
\Ascr=\ZAcH\ot_\ZAH A\quad{\rm and}\quad\Bscr=A_0\ot H^*
$$
as in the proof of Theorem 7.5. There is a homomorphism of $H$-module algebras
$$
\psi:\Ascr\to\Bscr,\qquad z\ot a\mapsto(z\ot1)\cdot\ph(a),
$$
where $\ph(a)=(\nu\ot\id)\bigl(\rho(a)\bigr)$. By Lemma 6.2 $\Ascr$ is 
$H$-simple. Hence $\Bscr$ is a free left $\Ascr$-module by Theorem 2.9. Let 
$r$ be its rank. Then
$$
P_{\Bscr/\ZAcH}\bigl(\ph(a),t\bigr)=P_{\Ascr/\ZAcH}(1\ot a,t)^r
=\nu^t P_{A/\ZAH}(a,t)^r
$$
for all $a\in A$. Here $\nu^t:\ZAH[t]\to\ZAcH[t]$ is the homomorphism of 
polynomial rings induced by the restriction of $\nu$ to $\ZAH$. Thus all 
coefficients of the polynomial above lie in $\nu\bigl(\ZAH\bigr)$.

Now let $z\in\ZAcH$. Pick any $a\in A$ such that $\nu(a)=z$. We have 
$H_0^*\cong H^*/J$ where $J$ is a Hopf ideal of $H^*$. The image of 
$\rho(a)\in A\ot H^*$ in $A_0\ot H_0^*$ coincides with $z\ot1$ since $z$ is 
$H_0$-invariant. In $A_0\ot H^*$ we get then $\ph(a)-z\ot1\in A_0\ot J$. But 
$J$ is nilpotent since $H_0$ contains the coradical of $H$. Hence 
$\ph(a)-z\ot1$ is nilpotent, and therefore
$$
P_{\Bscr/\ZAcH}\bigl(\ph(a),t\bigr)=P_{\Bscr/\ZAcH}(z\ot1,t)=(t-z)^m
$$
where $m=[\Bscr:\ZAcH]=[A_0:\ZAcH](\dim H)$. It follows that
$$
{m\choose j}z^j\in\nu\bigl(\ZAH\bigr)\quad\hbox{for all $j=0,1,\ldots,m$}.
$$
In particular, $mz\in\nu\bigl(\ZAH\bigr)$. If $\chr\bbk=0$, this entails 
$z\in\nu\bigl(\ZAH\bigr)$.

Suppose that $\chr\bbk=p>0$. Let $p^s$ be the largest power of $p$ dividing $m$. 
Taking $j=p^s$, we have ${m\choose j}\not\equiv0\mod p$, whence 
$z^j\in\nu\bigl(\ZAH\bigr)$.
\endproof

In \cite{Er-b16, Prop. 3.1} it was assumed that $H_0$ coincides with the 
coradical of $H$, but clearly the weaker assumption that $H_0$ contains the 
coradical of $H$ is sufficient.

We have given a slightly different proof. The proof in \cite{Er-b16} is based 
on the embedding of the simple $H$-module algebra $Z(A/M)\ot_\ZAH A$ into 
$A/M\ot H^*$ where $M$ is any maximal ideal of $A$ containing $K$. Using this 
embedding one can see that (ii) holds with $p^s$ taken to be the 
largest power of $p$ dividing the number
$$
[A/M:Z(A/M)]\cdot(\dim H).
$$
In other words, one obtains a possibly different value of $s$.

\proclaim
Corollary 8.2 \cite{Er-b16}.
Let $H_0\sbs H$ be a Hopf subalgebra containing the coradical of $H$. Suppose 
that $A$ is PI and prime {\rm(}or at least $H_0$-prime{\rm)}. If $\chr\bbk=0$, 
then
$$
Z(A)^{H_0}=\ZAH.
$$
\endproclaim

\Proof.
The quotient ring $Q=Q(A)$ is $H_0$-simple, and so Proposition 8.1 applies to 
the $H$-module algebra $Q$ and its ideal $K=0$.
\endproof

When $\chr\bbk>0$ Eryashkin has investigated the relationship between the 
central invariants for $H$ and for $H_0$ in $H$-prime PI algebras. The next 
result is a consequence of \cite{Er-b16, Prop. 3.2}:

\proclaim
Theorem 8.3.
Assume that $\chr\bbk=p>0$. Let $H_0\sbs H$ be a Hopf subalgebra containing the 
coradical of $H$. Suppose that $A$ is PI and $H$-prime. Let $A_0=A/P_0$ where 
$P_0$ is an $H_0$-prime ideal of $A$ containing no nonzero $H$-stable ideals 
of $A$. If $A_0$ is integral over $\ZAcH$, then $A$ is integral over $\ZAH$.
\endproclaim

In \cite{Er-b16, Th. 3.1} it was assumed additionally that $H_0$ is semisimple, 
which allows one to replace the condition that $A_0$ is integral over $\ZAcH$ 
with two weaker integrality assumptions for intermediate ring extensions.

\proclaim
Corollary 8.4 \cite{Er12}.
Assume that $\chr\bbk=p>0$ and $H$ is pointed. If $A$ contains a prime ideal $P$ 
such that the factor algebra $A/P$ is commutative and $P$ contains no nonzero 
$H$-stable ideals of $A$, then $A$ is integral over $\ZAH$.
\endproclaim

\Proof.
Here the coradical $H_0$ of $H$ is a group algebra $\bbk G$. With 
$P_0=\bigcap_{g\in G}gP$ we meet the hypothesis of Theorem 8.3 since 
$A_0=A/P_0$ is commutative, and so $A_0=Z(A_0)$ is integral over 
$A_0^{H_0}=A_0^G$ by the classical theory.
\endproof

\references
\nextref Art91
\auth{V.A.,Artamonov}
\paper{The structure of Hopf algebras\inRus}
\journal{Itogi Nauki Tekh. Ser. Algebra Topol. Geom.}
\Vol{29}
\Year{1991}
\Pages{3-63};
\etransl{J. Math. Sci.}
\Vol{71}
\Year{1991}
\Pages{2289-2328}

\nextref Art96
\auth{V.A.,Artamonov}
\paper{Invariants of Hopf algebras\inRus}
\journal{Vestnik Moscov. Univ. Ser. Mat. Mekh.}
\Vol{4}
\Year{1996}
\Pages{45-49};
\etransl{Moscow Univ. Math. Bull.}
\Vol{51}
\Year{1996}
\Pages{41-44}

\nextref Bah-L98
\auth{Yu.A.,Bahturin;V.,Linchenko}
\paper{Identities of algebras with actions of Hopf algebras}
\journal{J.~Algebra}
\Vol{202}
\Year{1998}
\Pages{634-654}

\nextref Ber-CF90
\auth{J.,Bergen;M.,Cohen;D.,Fischman}
\paper{Irreducible actions and faithful actions of Hopf algebras}
\journal{Isr. J. Math.}
\Vol{72}
\Year{1990}
\Pages{5-18}

\nextref Ber-M86
\auth{J.,Bergen;S.,Montgomery}
\paper{Smash products and outer derivations}
\journal{Isr. J. Math.}
\Vol{53}
\Year{1986}
\Pages{321-345}

\nextref Ber-I73
\auth{G.M.,Bergman;I.M.,Isaacs}
\paper{Rings with fixed-point-free group actions}
\journal{Proc. London Math. Soc.}
\Vol{27}
\Year{1973}
\Pages{69-87}

\nextref Bj71
\auth{J.-E.,Bj\"ork}
\paper{Conditions which imply that subrings of semiprimary rings are semiprimary}
\journal{J.~Algebra}
\Vol{19}
\Year{1971}
\Pages{384-395}

\nextref Ch90
\auth{W.,Chin}
\paper{Spectra of smash products}
\journal{Isr. J. Math.}
\Vol{72}
\Year{1990}
\Pages{84-98}

\nextref Coh86
\auth{M.,Cohen}
\paper{Smash products, inner actions and quotient rings}
\journal{Pacific J. Math.}
\Vol{125}
\Year{1986}
\Pages{45-66}

\nextref Coh92
\auth{M.,Cohen}
\paper{Hopf algebra actions -- revisited}
\journal{Contemp. Math.}
\Vol{134}
\Year{1992}
\Pages{1-18}

\nextref Coh-F86
\auth{M.,Cohen;D.,Fischman}
\paper{Hopf algebra actions}
\journal{J.~Algebra}
\Vol{100}
\Year{1986}
\Pages{363-379}

\nextref Coh-F92
\auth{M.,Cohen;D.,Fischman}
\paper{Semisimple extensions and elements of trace 1}
\journal{J.~Algebra}
\Vol{149}
\Year{1992}
\Pages{419-437}

\nextref Coh-FM90
\auth{M.,Cohen;D.,Fischman;S.,Montgomery}
\paper{Hopf Galois extensions, smash products, and Morita equivalence}
\journal{J.~Algebra}
\Vol{133}
\Year{1990}
\Pages{351-372}

\nextref Coh-R83
\auth{M.,Cohen;L.H.,Rowen}
\paper{Group graded rings}
\journal{Comm. Algebra}
\Vol{11}
\Year{1983}
\Pages{1253-1270}

\nextref Coh-W93
\auth{M.,Cohen;S.,Westreich}
\paper{Central invariants of $H$-module algebras}
\journal{Comm. Algebra}
\Vol{21}
\Year{1993}
\Pages{2859-2883}

\nextref Coh-W94
\auth{M.,Cohen;S.,Westreich}
\paper{From supersymmetry to quantum commutativity}
\journal{J.~Algebra}
\Vol{168}
\Year{1994}
\Pages{1-27}

\nextref Coh-WZ96
\auth{M.,Cohen;S.,Westreich;S.,Zhu}
\paper{Determinants, integrality and Noether's theorem for quantum commutative algebras}
\journal{Isr. J. Math.}
\Vol{96}
\Year{1996}
\Pages{185-222}

\nextref De-G
\auth{M.,Demazure;P.,Gabriel}
\book{Groupes Alg\'ebriques I}
\publisher{Masson}
\Year{1970}

\nextref Doi85
\auth{Y.,Doi}
\paper{Algebras with total integrals}
\journal{Comm. Algebra}
\Vol{13}
\Year{1985}
\Pages{2137-2159}

\nextref Er11
\auth{M.S.,Eryashkin}
\paper{Invariants of the action of a semisimple finite-dimensional Hopf algebra on special algebras\inRus}
\journal{Izv. Vyssh. Uchebn. Zaved. Mat.}
\Vol{8}
\Year{2011}
\Pages{14-22};
\etransl{Russian Math. (Iz. VUZ)}
\Vol{55}
\Year{2011}
\Pages{11-18}

\nextref Er12
\auth{M.S.,Eryashkin}
\paper{Martindale rings and $H$-module algebras with invariant characteristic polynomials\inRus}
\journal{Sibirsk. Mat. Zh.}
\Vol{53}
\Year{2012}
\Pages{822-838};
\etransl{Sib. Math. J.}
\Vol{53}
\Year{2012}
\Pages{659-671}

\nextref Er-a16
\auth{M.S.,Eryashkin}
\paper{Invariants and rings of quotients of $H$-semiprime $H$-module algebras satisfying a polynomial identity\inRus}
\journal{Izv. Vyssh. Uchebn. Zaved. Mat.}
\Vol{5}
\Year{2016}
\Pages{22-40};
\etransl{Russian Math. (Iz. VUZ)}
\Vol{60}
\Year{2016}
\Pages{18-34}

\nextref Er-b16
\auth{M.S.,Eryashkin}
\paper{Invariants of the action of a semisimple Hopf algebra on PI-algebra\inRus}
\journal{Izv. Vyssh. Uchebn. Zaved. Mat.}
\Vol{8}
\Year{2016}
\Pages{21-34};
\etransl{Russian Math. (Iz. VUZ)}
\Vol{60}
\Year{2016}
\Pages{17-28}

\nextref Et15
\auth{P.,Etingof}
\paper{Galois bimodules and integrality of PI comodule algebras over invariants}
\journal{J.~Noncommut. Geom.}
\Vol{9}
\Year{2015}
\Pages{567-602}

\nextref Et-W14
\auth{P.,Etingof;C.,Walton}
\paper{Semisimple Hopf actions on commutative domains}
\journal{Adv. Math.}
\Vol{251}
\Year{2014}
\Pages{47-61}

\nextref Et-W15
\auth{P.,Etingof;C.,Walton}
\paper{Pointed Hopf actions on fields. I}
\journal{Transform. Groups}
\Vol{20}
\Year{2015}
\Pages{985-1013}

\nextref Et-W16
\auth{P.,Etingof;C.,Walton}
\paper{Pointed Hopf actions on fields. II}
\journal{J.~Algebra}
\Vol{460}
\Year{2016}
\Pages{253-283}

\nextref Fer94
\auth{W.R.,Ferrer Santos}
\paper{Finite generation of the invariants of finite dimensional Hopf algebras}
\journal{J.~Algebra}
\Vol{165}
\Year{1994}
\Pages{543-549}

\nextref Kal-T08
\auth{M.,Kalniuk;A.,Tyc}
\paper{Geometrically reductive Hopf algebras and their invariants}
\journal{J.~Algebra}
\Vol{320}
\Year{2008}
\Pages{1344-1363}

\nextref Kato67
\auth{T.,Kato}
\paper{Self-injective rings}
\journal{Tohoku Math.~J.}
\Vol{19}
\Year{1967}
\Pages{485-495}

\nextref Kh74
\auth{V.K.,Kharchenko}
\paper{Galois extensions and rings of quotients\inRus}
\journal{Algebra i Logika}
\Vol{13}
\Year{1974}
\Pages{460-484};
\etransl{Algebra Logic}
\Vol{13}
\Year{1975}
\Pages{265-281}

\nextref Kr-T81
\auth{H.F.,Kreimer;M.,Takeuchi}
\paper{Hopf algebras and Galois extensions of an algebra}
\journal{Indiana Univ. Math.~J.}
\Vol{30}
\Year{1981}
\Pages{675-692}

\nextref Lin-M07
\auth{V.,Linchenko;S.,Montgomery}
\paper{Semiprime smash products and $H$-stable prime radicals for PI-algebras}
\journal{Proc. Amer. Math. Soc.}
\Vol{135}
\Year{2007}
\Pages{3091-3098}

\nextref Matc91
\auth{J.,Matczuk}
\paper{Centrally closed Hopf module algebras}
\journal{Comm. Algebra}
\Vol{19}
\Year{1991}
\Pages{1909-1918}

\nextref McC-R
\auth{J.C.,McConnell;J.C.,Robson}
\book{Noncommutative Noetherian Rings}
\publisher{Wiley}
\Year{1987}

\nextref Mo78
\auth{S.,Montgomery}
\paper{Outer automorphisms of semi-prime rings}
\journal{J.~London Math. Soc.}
\Vol{18}
\Year{1978}
\Pages{209-220}

\nextref Mo80
\auth{S.,Montgomery}
\book{Fixed Rings of Finite Automorphism Groups of Associative Rings}
\BkSer{Lecture Notes Math.}
\BkVol{818}
\publisher{Springer}
\Year{1980}

\nextref Mo93
\auth{S.,Montgomery}
\book{Hopf Algebras and their Actions on Rings}
\BkSer{CBMS Reg. Conf. Ser. Math.}
\BkVol{82}
\publisher{Amer. Math. Soc.}
\Year{1993}

\nextref Mo-S99
\auth{S.,Montgomery;H.-J.,Schneider}
\paper{Prime ideals in Hopf Galois extensions}
\journal{Isr. J. Math.}
\Vol{112}
\Year{1999}
\Pages{187-235}

\nextref Mo-Sm84
\auth{S.,Montgomery;L.W.,Small}
\paper{Integrality and prime ideals in fixed rings of P.I. rings}
\journal{J.~Pure Appl. Algebra}
\Vol{31}
\Year{1984}
\Pages{185-190}

\nextref Mum
\auth{D.,Mumford}
\book{Abelian Varieties}
\publisher{Oxford Univ. Press}
\Year{1970}

\nextref Par-S78
\auth{R.,Pare;W.,Schelter}
\paper{Finite extensions are integral}
\journal{J.~Algebra}
\Vol{53}
\Year{1978}
\Pages{477-479}

\nextref Pas
\auth{D.S.,Passman}
\book{Infinite Crossed Products}
\BkSer{Pure Appl. Math.}
\BkVol{135}
\publisher{Academic Press}
\Year{1989}

\nextref Qu89
\auth{D.,Quinn}
\paper{Integrality over fixed rings}
\journal{J.~London Math. Soc.}
\Vol{40}
\Year{1989}
\Pages{206-214}

\nextref Row
\auth{L.H.,Rowen}
\book{Ring Theory, Vol. II}
\publisher{Academic Press}
\Year{1988}

\nextref Sch76
\auth{W.,Schelter}
\paper{Integral extensions of rings satisfying a polynomial identity}
\journal{J.~Algebra}
\Vol{40}
\Year{1976}
\Pages{245-257};
Errata
\Vol{44}
\Year{1977}
576.

\nextref Sk02
\auth{S.,Skryabin}
\paper{Invariants of finite group schemes}
\journal{J.~London Math. Soc.}
\Vol{65}
\Year{2002}
\Pages{339-360}

\nextref Sk04
\auth{S.,Skryabin}
\paper{Invariants of finite Hopf algebras}
\journal{Adv. Math.}
\Vol{183}
\Year{2004}
\Pages{209-239}

\nextref Sk07
\auth{S.,Skryabin}
\paper{Projectivity and freeness over comodule algebras}
\journal{Trans. Amer. Math. Soc.}
\Vol{359}
\Year{2007}
\Pages{2597-2623}

\nextref Sk-ART11
\auth{S.,Skryabin}
\paper{Structure of $H$-semiprime Artinian algebras}
\journal{Algebr. Represent. Theory}
\Vol{14}
\Year{2011}
\Pages{803-822}

\nextref Sk-JA11
\auth{S.,Skryabin}
\paper{Coring stabilizers for a Hopf algebra coaction}
\journal{J.~Algebra}
\Vol{338}
\Year{2011}
\Pages{71-91}

\nextref Sk15
\auth{S.,Skryabin}
\paper{Invariant subrings and Jacobson radicals of Noetherian Hopf module algebras}
\journal{Isr. J. Math.}
\Vol{207}
\Year{2015}
\Pages{881-898}

\nextref Sk-PAMS17
\auth{S.,Skryabin}
\paper{Finiteness of the number of coideal subalgebras}
\journal{Proc. Amer. Math. Soc.}
\Vol{145}
\Year{2017}
\Pages{2859-2869}

\nextref Sk-JA17
\auth{S.,Skryabin}
\paper{The left and right dimensions of a skew field over the subfield of invariants}
\journal{J.~Algebra}
\Vol{482}
\Year{2017}
\Pages{248-263}

\nextref Sk-Oy06
\auth{S.,Skryabin;F.,Van Oystaeyen}
\paper{The Goldie theorem for $H$-semiprime algebras}
\journal{J.~Algebra}
\Vol{305}
\Year{2006}
\Pages{292-320}

\nextref Tot97
\auth{A.A.,Totok}
\paper{On invariants of finite-dimensional pointed Hopf algebras\inRus}
\journal{Vestnik Moscov. Univ. Ser. Mat. Mekh.}
\Vol{3}
\Year{1997}
\Pages{31-34};
\etransl{Moscow Univ. Math. Bull.}
\Vol{52}
\Year{1997}
\Pages{33-36}

\nextref Tot98
\auth{A.A.,Totok}
\paper{Actions of Hopf algebras\inRus}
\journal{Mat. Sbornik}
\Vol{189}
\Year{1998}
\Pages{149-160};
\etransl{Sb. Math.}
\Vol{189}
\Year{1998}
\Pages{149-159}

\nextref Zhu96
\auth{S.,Zhu}
\paper{Integrality of module algebras over its invariants}
\journal{J.~Algebra}
\Vol{180}
\Year{1996}
\Pages{187-205}

\endreferences
\bye